\documentclass[12pt]{article}

\usepackage[]{pdfsync} 
\usepackage[francais,english]{babel} 
\usepackage[]{latexsym} 
\usepackage{amsmath} 
\usepackage{graphicx} 
\usepackage{amsfonts} 
\usepackage{amssymb}
\usepackage[T1]{fontenc}
\if{ 
\usepackage{fourier} 
\usepackage[scaled=0.875]{helvet} 

\usepackage{lmodern} 
\usepackage[T1]{fontenc}

\usepackage{mathpazo} 
\linespread{1.05}        
\usepackage[scaled]{helvet} 
\usepackage{courier} 
\normalfont
\usepackage[T1]{fontenc}

\linespread{1.05}        
\usepackage[scaled]{helvet} 
\usepackage{courier} 
\usepackage{euler} 
\normalfont
\usepackage[T1]{fontenc}
} \fi 

\newcommand \as{almost surely }  

\font \db = msbm10 at 12 pt
\font \sdb = msbm7 at 8.5 pt

      \def \t{\theta}  
\def \e{\varepsilon}  \def \f{\varphi}    \def \rr{\varrho}

\font \db = msbm10 at 12 pt
\font \sdb = msbm7 at 8.5 pt     \font \bdb = msbm10 at 14 pt  
\def \R{\mbox{\db R}}  \def \H{\mbox{\db H}} \def \E{\mbox{\db E}} \def \P{\mbox{\db P}} 
\def \N{\mbox{\db N}} 
\def \sR{\mbox{\sdb R}}    
\def \sN{\mbox{\sdb N}} 
\def \bR{\mbox{\bdb R}}     \def \bH{\mbox{\bdb H}}

\newcommand \parsn{\par \smallskip \noindent }  \newcommand \parmn{\par \medskip \noindent } 
\newcommand \pars{\par \smallskip }  \newcommand \parm{\par \medskip }

   \def \C{\mbox{\db C}}   \def \Z{\mbox{\db Z}}

\def \FF{{\cal F}}  \def \GG{{\cal G}}         
\def \LL{{\cal L}}  \def \M{{\cal M}}       \def \O{{\cal O}}    \def \PP{{\cal P}}
  \def \RR{{\cal R}}

  \def \lra{\longrightarrow}   \def \ub{\underbar}        \def\lmt{\longmapsto}
\def \Ra{\Rightarrow}     \def \LRa{\Leftrightarrow} 
\def\LRA{\Longleftrightarrow}
              \def \nea{\nearrow}        \def \sea{\searrow}

   \def \ds{\displaystyle}   \def \ts{\textstyle}       \def \ss{\scriptstyle}
\def \rt1{\sqrt{-1}\,\,}  \def \1{^{-1}}            \def \2{^{-2}}             \def \5{{\ts \frac{1}{2}}}
\def\moins{\rotatebox{30}{\raise0.3pt\hbox{$\,\scriptstyle\setminus$}\,}} 

\def\cotg{{\rm cotg}\,}   \def\tg{{\rm tg}\,}   \def\coth{{\rm coth}\,}        \def\th{{\rm th}\,}
\def\ch{{\rm ch}\,}       \def\sh{{\rm sh}\,}      

\def \parn{\par\noindent }  \def\indf{\leavevmode\indent }  
\def \parsn{\par\smallskip \noindent } \def \parmn{\par\medskip \noindent } 
\def \parm{\par \medskip } \def \pars{\par \smallskip }

\begin{document}

\newtheorem{defi}[subsubsection]{Definition} 
\newtheorem{theo}[subsubsection]{Theorem} 
\newtheorem{prop}[subsubsection]{Proposition}  
\newtheorem{cor}[subsubsection]{Corollary} 
\newtheorem{lem}[subsubsection]{Lemma} 
\newtheorem{rem}[subsubsection]{Remark} 
\newtheorem{subs}[subsubsection]{} 

\newtheorem{Cor}[subsection]{Corollary} 
\newtheorem{Lem}[subsection]{Lemma} 
\newtheorem{Theo}[subsection]{Theorem} 
\newtheorem{Prop}[subsection]{Proposition} 
\newtheorem{Rem}[subsection]{Remark} 

\newcommand \beq{\begin{equation}} \newcommand \eeq{\end{equation}} 

\newcommand \bthe{\begin{theo}}  \newcommand \ethe{\end{theo}}   
\newcommand \bpro{\begin{prop}}  \newcommand \epro{\end{prop}}   
\newcommand \bcor{\begin{cor}}    \newcommand \ecor{\end{cor}}     
\newcommand \blem{\begin{lem}}   \newcommand \elem{\end{lem}}   
\newcommand \brem{\begin{rem}}  \newcommand \erem{\end{rem}} 
\newcommand \bdefi{\begin{defi}}   \newcommand \edefi{\end{defi}}  
\newcommand \bsub{\begin{subs} \  }   \newcommand \esub{\end{subs}}  

\newcommand \bCor{\begin{Cor}}    \newcommand \eCor{\end{Cor}}     
\newcommand \bLem{\begin{Lem}}   \newcommand \eLem{\end{Lem}}   
\newcommand \bThe{\begin{Theo}}  \newcommand \eThe{\end{Theo}}   
\newcommand \bPro{\begin{Prop}}  \newcommand \ePro{\end{Prop}}   
\newcommand \bRem{\begin{Rem}}  \newcommand \eRem{\end{Rem}}

\title{\bf Small time asymptotics for an example of strictly hypoelliptic heat kernel}
\author{Jacques FRANCHI }
\date {August 2012}
\maketitle

\begin{abstract}
A small time asymptotics of the density is established for a simplified (non-Gaussian, strictly hypoelliptic) second chaos process tangent to the Dudley relativistic  diffusion.
\end{abstract} 

\section{Introduction} \indf 
   The problem of estimating the heat kernel, or the density of  a diffusion, particularly as time goes to zero, has been extensively studied for a long time. Let us mention only the articles [V], [A], [BA1], [BA2], [L], [ERS], and the existence of other works on that subject by Azencott, Molchanov and Bismut, quoted in  [BA1]. \par
   To summary roughly, a very classical question addresses the asymptotic behavior (as $\,s\sea 0$) of the density $\,p_s(x,y)$ of the diffusion $(x_s)$ solving a Stratonovich stochastic differential equation 
\begin{equation*}
x_s = x + \sum_{j=1}^k \int_0^s V_j(x_\tau)\circ dW^j_\tau + \int_0^s V_0(x_\tau)\, d\tau \, , 
\end{equation*}
where the smooth vector fields $\,V_j\,$ are supposed to satisfy a H\"ormander condition\;; the underlying space being $\R^d$ or some $d$-dimensional smooth manifold $\M$. \par
   The elliptic case being very well understood for a long time ([V], [A]), the studies focussed then on the sub-elliptic case, that is to say, when the strong H\"ormander condition (that the Lie algebra generated by the fields $\,V_1,\ldots, V_k\,$ has maximal rank everywhere) is fulfilled. In that case these fields generate a sub-Riemannian distance $\,d(x,y)$, defined as in control theory, by considering only $C^1$ paths 
whose tangent vectors are spanned by them. Then the wanted asymptotic expansion tends to have the following Gaussian-like form\,: 
\begin{equation} \label{f.GaussLAsymp} 
p_s(x,y) = s^{-d/2}\, \exp\!\big( - d(x,y)^2/(2 s)\big) \bigg( \sum_{\ell=0}^n \gamma_\ell(x,y)\, s^\ell + \O(s^{n+1})\bigg)  
\end{equation}
for any $\,n\in\N^*$, with smooth $\,\gamma_\ell$'s and $\gamma_0>0\,$, provided $(x,y)$ does not belong to the cut-locus (and uniformly within any compact set which does not intersect the cut-locus). See in particular ([BA1], th\'eor\`eme 3.1). Note that the condition of remaining outside the cut-locus is here necessary, as showed in particular by [BA2]. \par
   The methods used to get this or a similar result have been of different nature. In [BA1], G. Ben Arous proceeds by expanding the flow associated to the diffusion (in this direction, see also [Ca]
) and using a Laplace method applied to the Fourier transform of $\,x_s\,$, then inverted by means of Malliavin's calculus. \par 
   The strictly hypoelliptic case, i.e., when only the weak H\"ormander condition (requiring the use of the drift vector field $V_0$ to recover the full tangent space) is fulfilled, remains much more problematic, and then rarely addressed. There is a priori no longer any reason that in such case the asymptotic behavior of $\,p_s(x,y)$ remains of the Gaussian-like type (\ref{f.GaussLAsymp}), all the less as a natural candidate for replacing the sub-Riemannian distance $\,d(x,y)$ is missing. Indeed this already fails for the mere (however Gaussian) Langevin process $\big(\beta_s,\int_0^s\beta_\tau d\tau\big)$\,: the missing distance is replaced by a time-dependent distance which presents some degeneracy in one direction, namely $\,\frac{6}{s^3}\big|(x-y) - \frac{s}{2}\, (\dot x-\dot y)\big|^2 + \frac{1}{2s}|\dot x-\dot y|^2\,$\raise0pt\hbox{,} see (\ref{f.Mehl}), (\ref{f.qtildes}) below.  See also [DM] for a more involved (non-curved, strictly hypoelliptic, perturbed) case where Langevin-like estimates hold (without precise asymptotics), roughly having the following Li-Yau-like form\,:
\begin{equation} \label{F.DM} 
C\1\, s^{-N}\,e^{-C\, d_s(x_s,y)^2} \,\le \, p_s(x,y)\, \le \, C\, s^{-N}\,e^{-C\1 \, d_s(x_s,y)^2} \,, \quad \hbox{for }\; 0<s<s_0\,. 
\end{equation}
    \pars
   
   In this article, an interesting case of rather natural hypoelliptic diffusion is considered first\,: that of a relativistic diffusion, first constructed in Minkowski's space (see [Du]), which makes sense on a generic smooth Lorentzian manifold, see [F-LJ1], [F-LJ2], [F-LJ3]. In the simplest case of Minkowski's space, it consists in the pair $(\xi_s, \dot\xi_s)\in \R^{1,d}\times\H^d$ (parametrized by its proper time $s$, and analogous to a Langevin process), where the velocity $(\dot\xi_s)$ is a hyperbolic Brownian motion. In the general case this Dudley diffusion can be rolled without slipping from a reference tangent space to the Lorentzian manifold, see [F-LJ1]. Note that even in the Minkowski space, there is a curvature constraint to be taken into account, namely that of the mass shell $\,\H^d$, at the heart of this framework. Moreover the relativistic diffusion is never sub-elliptic, but only hypoelliptic, and a priori a Gaussian-like asymptotic expansion as (\ref{f.GaussLAsymp}) does not even make sense, since there is no longer any natural candidate to replace the sub-Riemannian distance $\,d(x,y)$. See however [BF], where some non-trivial information is extracted about the relativistic diffusion, by considering the sub-Riemannian distance generated by all fields $V_0,V_1,\ldots, V_k$ (i.e., not only $V_1,\ldots, V_k$). Talking of this, an important feature of the strictly hypoelliptic case, which is fulfilled in the relativistic framework, is when the graded geometry generated by the successive brackets of a given weight is (at least locally) constant (see [NSW], and also [T]), yielding homogeneity in the afore-mentioned time-dependent distance. 
\pars 
 
  To proceed, we shall compute a Fourier-Laplace transform, which  seems to be the only way of getting any quantitative access to the density kernel of the relativistic diffusion \big([BA1] relies already on the Fourier transform, but then the method followed by G. Ben Arous is based on a stochastic variation about a minimal geodesic, which does not exist here, and on the local strict convexity of the energy functional (due to the sub-ellipticity), which does not hold here\big). Because of the singularity in the most natural polar coordinates, we first choose alternative, less intuitive but smooth coordinates, and using them, partially expand the relativistic diffusion to project it on the second Wiener chaos, thereby exhibiting a simplified ``tangent process''. For this simplified process the Fourier-Laplace transform is exactly computable. Then analyzing its inverse Fourier transform very carefully and using a saddle-point method allows to derive an asymptotic equivalent for the density $\,q_s\,$ of this tangent process, as time $\,s\,$ goes to zero (see Theorem \ref{pro.asymptqs} below). As in the Langevin or in the more sophisticated case of [DM], the exponential term is given by a time-dependent distance, namely the same as the afore-mentioned Langevin one, the strictly second chaos  coordinate appearing only in the off-exponent term. The initial analogous question about the relativistic diffusion remains open, as the degree of contact between both considered processes (the effective computation of the Fourier-Laplace  transform being bounded to the second Wiener chaos) seems so far too weak, to allow to deduce a former asymptotic behaviour from the second one. The non-appearance of the non-Gaussian coordinates in the found asymptotic exponent lets however think that this could remain so for all higher order chaos terms of the Taylor expansion. Thence a tempting guess, resulting from both the sub-elliptic case (\ref{f.GaussLAsymp}) as solved by [BA1], the modified Li-Yau-like estimates (\ref{F.DM}) obtained in [DM] and the present work, whose main result is Theorem \ref{pro.asymptqs} below, is that an expansion having the following form could (maybe generally, under consistency of the Lie graded geometry) hold\,:
\begin{equation*}
p_s(x,y) =\,s^{-N}\, \exp\!\Big( - d_s(x,y)^2 \Big) \bigg( \sum_{\ell=0}^n \gamma_\ell(x,y)\, s^\ell + \O(s^{n+1})\bigg) . 
\end{equation*}
 \par

   In this first attempt we restrict to the simplest case of the five-dimensional Minkowski space $\,\R^{1,2}\times\H^2$ and then to its five-dimensional second chaos tangent process. The case of the generic Minkowski space $\,\R^{1,d}\times\H^d$ is actually very analogous, but would mainly bring notational difficulties without modifying the method.  We hope that this particular toy example will allow to understand better what can happen and could be undertaken, concerning small time asymptotics of the relativistic diffusion itself, and then maybe in some more generic strictly hypoelliptic framework. \par 
   The content is organized as follows. \parn
In Section \ref{sec.RepHR2} are mainly described the setting and the smooth parametrization used then.    
In Section \ref{sec.sta} the simplified ``tangent process'' $(Y_s)$ to the relativistic diffusion $(X_s)$ is exhibited.  
In Section \ref{sec.FLTransf} the Fourier-Laplace transform of the tangent process $(Y_s)$ is computed. 
In Section \ref{s.intexprdqs} a closed integral expression for the density $q_s$ of $(Y_s)$ is explicited.
Section \ref{s.asbdqs} yields a precise (off the diagonal) equivalent for the density $\,q_s\,$ as $\,s\to 0$, which is the content of the main result (Theorem \ref{pro.asymptqs}).   \   
Section \ref{sec.Proofs} contains three rather technical proofs, which have been postponed till there to lighten the reading. 


\section{A smooth parametrization of $\,\bH^2\times\bR^{1,2}\equiv T^1_+\bR^{1,2}$} \label{sec.RepHR2} \indf 
   The Dudley relativistic diffusion $X_s = (\dot\xi_s,\xi_s)$ (see [Du]) lives in the future-directed unit tangent bundle $T^1_+\R^{1,2}\equiv \H^2\times\R^{1,2}$ to the Minkowski space $\big(\R^{1,2}, \langle \cdot , \cdot\rangle\big)$ \big(or alternatively, in its frame bundle, isomorphic to the Poincar\'e isometry group $\,\PP^3= \mathrm{PSO}(1,2)\propto\R^{1,2}$\big). We classically identify the hyperbolic plane $\H^2\subset \R^{1,2}$ with the upper sheet of the hyperboloid having equation $\,|\dot\xi^0|^2-|\dot\xi^1|^2-|\dot\xi^2|^2 \equiv \langle \dot\xi, \dot\xi\rangle =1$ within  $\R^{1,2}$ \big(endowed with its canonical basis $(e_0,e_1,e_2)$\big). The velocity sub-diffusion $(\dot\xi_s)$ is a hyperbolic Brownian motion, and we merely have $\,d\xi_s = \dot\xi_s\, ds\,$. The parameter $\,s\,$ is precisely the physical proper time. \pars 

   We shall use the following smooth  coordinates $(\lambda\,,\, \mu\,,\, x\,,\, y\,,\, z)\in\R^5$ on $T^1_+\R^{1,2}\,$: 
\begin{equation}  \label{f.nouvcoord} 
\dot\xi^0 = \ch \lambda\,\ch\mu \; ;\quad \dot\xi^1 = \ch \lambda\,\sh\mu \; ;\quad \dot\xi^2 = \sh \lambda \; ;\quad  \xi^0 = x \; ;\quad \xi^1 = y\; ;\quad  \xi^2 = z\, . 
\end{equation}
In these coordinates the Dudley diffusion $\,X_s \equiv (\lambda_s\,,\, \mu_s\,,\, x_s\,,\, y_s\,,\, z_s)\in\R^5$ satisfies the following system of stochastic differential equations (for independent real Brownian motions $\,w,\beta$)$\,$: 
\begin{equation}  \label{f.sde} 
d\lambda_s = \sigma\, dw_s + {\ts\frac{\sigma^2}{2}}\, \th\lambda_s \, ds \; ;\quad d\mu_s =  \sigma\, \frac{d\beta_s}{\ch\lambda_s}\; ; 
\end{equation}
\begin{equation}  \label{f.sde'} 
dx_s = \ch \lambda_s\,\ch\mu_s\, ds  \; ;\quad dy_s = \ch \lambda_s\,\sh\mu_s\,ds \; ;\quad dz_s = \sh \lambda_s\, ds\, . 
\end{equation}
The infinitesimal generator $\,\LL\,$ of $(X_s)$ reads in these coordinates$\,$: 
\begin{equation}  \label{f.gen} 
\LL =\, {\ts\frac{\sigma^2}{2}} \bigg[ \frac{\partial^2}{\partial\lambda^2} + \th\lambda\, \frac{\partial}{\partial\lambda} + (\ch\lambda)\2 \frac{\partial^2}{\partial\mu^2} \bigg] + \ch \lambda\,\ch\mu \, \frac{\partial}{\partial x} + \ch \lambda\,\sh\mu \, \frac{\partial}{\partial y} + \sh \lambda\,\frac{\partial}{\partial z}\, \raise1.9pt\hbox{.} 
\end{equation}
Consider the following smooth vector fields on $\R^5$: 
$$V_1:= \sigma\,\frac{\partial}{\partial\lambda}\, \raise1.8pt\hbox{,}  \,V_2:= \frac{\sigma}{\ch \lambda}\,\frac{\partial}{\partial\mu}\, \raise1.8pt\hbox{,}  \   V_0' := \frac{\partial}{\partial s} + {\ts\frac{\sigma^2}{2}}\,\th\lambda\,\frac{\partial}{\partial\lambda} + \ch\lambda\,\ch\mu\, \frac{\partial}{\partial x} + \ch\lambda\,\sh\mu\, \frac{\partial}{\partial y} + \sh\lambda\, \frac{\partial}{\partial z}\, \raise1.8pt\hbox{,} $$ 
and also \quad $V_0:= V'_0 - \frac{\partial}{\partial s} \,$ \raise1.8pt\hbox{.} \quad We have then  
\pars 
\centerline{$ \frac{\partial}{\partial s} + \LL =\,\5\big( V_1^2+V_2^2\big) + V'_0\,, \ [V_1,V'_0] = \frac{\sigma^3}{2\,\ch^2\lambda}\,\frac{\partial}{\partial\lambda} + \sigma\,\sh\lambda\,\ch\mu\, \frac{\partial}{\partial x} + \sigma\,\sh\lambda\,\sh\mu\, \frac{\partial}{\partial y} + \sigma\,\ch\lambda\, \frac{\partial}{\partial z}\,\raise1.5pt\hbox{,} 
$}\pars 
\centerline{$ [V_2\,,V'_0] = \frac{\sigma^3\, \th^2\lambda}{2\,\ch\lambda}\,\frac{\partial}{\partial\lambda}  + \sigma\,\sh\mu\, \frac{\partial}{\partial x} + \sigma\,\ch\mu\, \frac{\partial}{\partial y}  \,\raise1.5pt\hbox{,} 
\quad \big[V_2 , [V_2,V'_0]\big] =  \frac{\sigma^2}{\ch\lambda}\big(\ch\mu\, \frac{\partial}{\partial x} +\sh\mu\, \frac{\partial}{\partial y}\big) \raise1.1pt\hbox{,} 
$}\parsn  
and then $\Big( V'_0,V_1, V_2, [V_1,V'_0] ,  [V_2,V'_0] , \big[V_2 , [V_2,V'_0]\big]\Big)$ has full rank 6 at any point, so that the weak H\"ormander condition holds. 
\parn 
Hence by hypoellipticity, $\,\LL\,$ admits a smooth heat kernel $\,p_s(X_0;X)$  ($s\in\R_+\,$,  $\,X_0,X\in\R^5$), with respect to the Liouville measure $\,L\,$, which reads \   $L(dX) = \ch\lambda\, d\lambda\, d\mu \, dx\,dy\,dz\,$: 
$$\E_{X_0}\big[F(X_s)\big] = \int_{T^1_+\sR^{1,2}} p_s(X_0;X)\, L(dX) = \int_{\sR^5} p_s(X_0;\,\lambda \,,\, \mu \,,\, x \,,\, y \,,\, z )\, \ch\lambda\, d\lambda\, d\mu \, dx\,dy\,dz\, . $$  
\par 

   An open question is to estimate \  $p_s(X_0\,;X)$, for small proper times $\,s\,$.\parsn
Up to apply some element of the Poincar\'e group $\,\PP^3$, we can restrict to $\,X_0 = (e_0,0) \equiv (0,0,0,0,0)$. Thus we have to deal with \  $p_s(0,X)\equiv p_s(X)$, \  for \  $X \equiv (\lambda \,,\, \mu \,,\, x \,,\, y \,,\, z )\in\R^5$. \parm 


   The underlying unperturbed (deterministic) process $X^0 \equiv (\lambda^0 \,,\, \mu^0 \,,\, x^0 \,,\, y^0 \,,\, z^0 )\in\R^5$  solves$\,$: 
\begin{equation*}  \label{f.sde0} 
d\lambda^0_s = {\ts\frac{\sigma^2}{2}}\, \th\lambda^0_s \, ds \, ;\quad \mu^0_s = 0\, ; \quad dx^0_s = \ch \lambda^0_s\, ds\,;\quad y^0_s = 0\, ;\quad dz^0_s = \sh \lambda^0_s\, ds\,  , 
\end{equation*}
and then is merely given by the geodesic \   $X^0_s \equiv (0 \,,\, 0 \,,\, s \,,\, 0 \,,\, 0 )$ (for any proper time $s$). \parsn  
Up to change the speed of the canonical Brownian motion $(w,\beta)$, by considering $(w_{\sigma^2 s},\beta_{\sigma^2 s})$ instead of $(w_s,\beta_s)$, we can absorb the speed parameter $\,\sigma\,$, and then suppose that $\,\sigma=1\,$.  \par 

\section{A process tangent to the relativistic diffusion $(X_s)$} \label{sec.sta} \indf 
    The main Theorem 2.1 in [Ca] could apply here (beware however that $V_0$ is unbounded), yielding a full general Taylor expansion for the diffusion $(X_s)$, in terms of the above vector fields $V_1,V_2,V_0$, their successive brackets, and of the iterated Stratonovich integrals with respect to $(w,\beta)$. Indeed, Equations (\ref{f.sde}),(\ref{f.sde'}) read equivalently\,:  
\begin{equation*}
dX_s = V_1(X_s)\, d w_s + V_2(X_s)\, d\beta_s + V_0'(X_s)\, ds\,. 
\end{equation*}
In [Ca], the successive remainders corresponding to the truncated Taylor expansion are controlled in probability. In this spirit and also almost surely, the process $X_s$ is approached as follows. 
\bLem  \label{lem.approxXY} \  $(i)$ \  For any $\,\e>0$, almost surely as proper time $\,s\sea 0\,$ we have\,: 
\begin{equation*}
X_s = \Big( w_s + {\ts\frac{1}{2}}\!\int_0^s\! w_\tau\,d\tau , \; \beta_s + o(s^{3/2-\e}) , \; s+ {\ts\frac{1}{2}}\! \int_0^s\! \big[\beta_\tau^2 + w_\tau^2\big] d\tau , \int_0^s\! \beta_\tau\,d\tau , \int_0^s\! w_\tau\,d\tau \Big) + o(s^{5/2-\e}) . 
\end{equation*}
$(ii)$ \  Setting \quad  ${\ds R'_s := (x_s , y_s , z_s) - \Big( s+ {\ts\frac{1}{2}}\! \int_0^s\! \big[\beta_\tau^2 + w_\tau^2\big] d\tau , \int_0^s\! \beta_\tau\,d\tau , \int_0^s\! w_\tau\,d\tau \Big) }$ \  and \parn
 ${\ds R_s := (\lambda_s , \mu_s) - ( w_s ,\beta_s)}$, \ 
there exist $\,c,\kappa >0\,$ such that \   for any $\,R>c\,$ we have\,:   \parn   
\centerline{${\ds \lim_{s\sea 0}\,\P\bigg[ \sup_{0\le t\le s}\, \| R_t\| \ge R\, s^{3/2}\bigg] \le\, e^{-R^\kappa/c}  }$ \   and \quad ${\ds \lim_{s\sea 0}\,\P\bigg[ \sup_{0\le t\le s}\, \| R'_t\| \ge R\, s^{5/2}\bigg] \le\, e^{-R^\kappa/c}}$.} 
\eLem 
The proof is posponed to Section \ref{sec.Proofs}.  \par
\bRem \label{rem.approxmus} \  {\rm  More precisely, concerning the martingale $(\mu_s)$ we have 
\begin{equation*}  
\mu_s =  \int_0^s \frac{d\beta_\tau}{\ch\lambda_\tau} = \int_0^s \big(1-\5\, w_\tau^2 + o(\tau^{2-\e})\big) d\beta_\tau\, = \beta_s -\5 \int_0^s w_\tau^2\, d\beta_\tau + o(s^{5/2-\e}) .
\end{equation*}
But the method used then does not work with the non-quadratic martingale ${\ds \int_0^s w_\tau^2\, d\beta_\tau}$\par
\centerline{\big(which equals $\, w_s^2 \beta_s - 2 \int_0^s \beta_\tau w_\tau dw_\tau -  \int_0^s \beta_\tau d\tau = o(s^{3/2-\e}) \big).$}
} \eRem
\parm  

 As a consequence, we shall use a perturbation method, approaching (for small proper time $s$) the relativistic diffusion $X_s\,$ by means of the $\R^5$-valued ``tangent process''\,:   
\begin{equation}  \label{f.2chaos} 
Y_s := \Big( w_s\, , \; \beta_s\, , \; {\ts\frac{1}{2}} \int_0^s \big[\beta_\tau^2 + w_\tau^2\big] d\tau\, , \, \int_0^s \beta_\tau\,d\tau\, , \, \int_0^s w_\tau\,d\tau \Big) =: \big( w_s\, , \, \beta_s\, , \,A_s\, , \,\zeta_s\, , \,\bar z_s \big)  
\end{equation}
which is not Gaussian, but has its third coordinate $\,A_s\,$ in the second Wiener chaos.  This actually yields the orthogonal projection of the process $(X_s)$ onto the second Wiener chaos. 
\bRem  \label{rem.GaussDM} \  {\rm  The fact that the second chaos term is needed in the approximation (without it, the tangent process would clearly not admit any density) makes a significant difference with the situation exhaustively investigated in [DM], where the approaching process is Gaussian. This can no longer be the case in the present setting, though both settings share the feature of being strictly hypoelliptic. A difference between both is the curvature, at the heart of the relativistic realm (even in the present Minkowski-Dudley flat case), due to the mass shell constraint on velocities. 
}\eRem
Note that for any fixed proper time $\,s>0\,$ we have$\,$: 
\begin{equation}  \label{f.2chaos'1} 
Y_{s}  \,\stackrel{law}{\equiv} \,\Big(  \sqrt{s}\, w_1\, , \;  \sqrt{s}\, \beta_1\, , \; {\ts\frac{s^2}{2}} \int_0^1 \big[\beta_\tau^2 + w_\tau^2\big] d\tau\, , \;  \sqrt{s^3}  \int_0^1 \beta_\tau\,d\tau\, , \; \sqrt{s^3} \int_0^1 w_\tau\,d\tau \Big) . 
\end{equation}
Denote by \  $ q_s= q_s(w,\beta , x , \zeta , z)\,$ the density of $Y_s\,$ with respect to the Lebesgue measure $\,\Lambda(dY) = dw\, d\beta \, dx\,d\zeta\,dz$ (i.e., {\it not} the Liouville measure $\,L$) on $\,\R^2\times\R_+^*\times\R^2$. \parn 
By the scaling property (\ref{f.2chaos'1}), it must satisfy : 
\begin{equation} \label{f.densYs} 
q_s(w,\beta , x , \zeta , z) = \frac{1}{s^6}\,\, q_1\bigg( \frac{w}{\sqrt{s}}\, \raise1.8pt\hbox{,}\, \frac{\beta}{\sqrt{s}}\, \raise1.8pt\hbox{,}\,  \frac{x}{s^2}\, \raise1.8pt\hbox{,}\, \frac{\zeta}{\sqrt{s^3}}\, \raise1.8pt\hbox{,}\, \frac{z}{\sqrt{s^3}} \bigg) \, , 
\end{equation}
and otherwise$\,$: \quad  
$ q_s(w,\beta , x , \zeta , z) =  q_s(-w,\beta , x , \zeta, -z) = q_s(w,-\beta , x , -\zeta , z)\, $. \par 


\if{ 
\blem  \label{lem.decrq1} \  The density $\,q_1= q_1(w,\beta , x , \zeta , z)$ of the variable $Y_1\,$ decreases eventually with respect to $\,x\,$:  \quad \parn 
there exists some positive function $\,\f = \f(w,\beta , \zeta , z)= \f(-w,\beta , \zeta , -z)= \f(w,-\beta , -\zeta , z)$ on $\R^4$ such that $\,x\mapsto q_1(w,\beta , x , \zeta , z)$ decreases on $\big[\f(w,\beta , \zeta , z),\infty\big[\,$. 
??? \  Actually  
\begin{equation*}
\f(w,\beta , \zeta , z)\, = \, \E\big[ A_1\,\big|\, w_1=w ,\,\beta_1=\beta,\,\zeta_1=\zeta,\,\bar z_1=z\big] \hbox{ (computed in Remark \ref{rem.Psi'0} below)}. 
\end{equation*}  
\elem 
\ub{Proof} \quad ?? 
}\fi 

\section{Fourier-Laplace transform of the tangent process} \label{sec.FLTransf}  
\subsection{Fourier-Laplace transform of the simplified process $Z_s$} \label{sec.DensZs} \indf 
    We need information on the density at time $s\,$ of the tangent process $(Y_s)$ given in (\ref{f.2chaos}). By the independence of $\,w,\beta\,$, it will be enough to consider the density of 
\begin{equation}  \label{f.2chaoss} 
Z_s := \Big( w_s\, , \, \int_0^s w_\tau\,d\tau\, , \, \int_0^s w_\tau^2\, d\tau\Big) \,\stackrel{law}{\equiv} \, \Big(  \sqrt{s}\, w_1\, , \, \sqrt{s^3} \int_0^1 w_\tau\,d\tau\, , \, {s^2} \int_0^1 w_\tau^2\, d\tau\Big) . 
\end{equation}
Note that this simplified tangent process $\,Z_s\,$ does not have any component beyond the second chaos.
     Because of the scaling property (\ref{f.2chaoss}) of  $\,Z_s\,$, its density $\,q^0_s(w,z,x)$ satisfies 
\begin{equation*}  
q^0_s(w, z, x) = \frac{1}{s^4}\,\, q^0_1\Big( \frac{w}{\sqrt{s}}\, \raise1.8pt\hbox{,}\, \frac{z}{\sqrt{s^3}}\, \raise1.8pt\hbox{,}\,  \frac{x}{s^2}\Big)\; =   q^0_s(-w, -z, x) \, . 
\end{equation*}
   Of course, the Langevin process ${\ds \Big( w_s\, ,  \int_0^s w_\tau\,d\tau\Big)}$ is Gaussian with covariance \parn
${ K^0_{s} = \begin{pmatrix} s & s^2/2 \cr s^2/2 &  s^3/3  \end{pmatrix}}$, so that it has the well-known  density 
\begin{equation} \label{f.Mehl} 
(w,z)\lmt \,\frac{\sqrt{3}}{\pi\, s^2}\, e^{-(6z^2-6s\, zw + 2s^2 w^2)/s^3} = \frac{\sqrt{3}}{\pi\, s^2}\, \exp\!\bigg[ - \frac{6}{s^3}\Big(z - \frac{s}{2}\, w\Big)^2 - \frac{w^2}{2s}\bigg] . 
\end{equation}
\centerline{\big(In particular, the expected value of $\int_0^s w_\tau\,d\tau$, conditionally on $\,w_s=w$, equals $\,sw/2$.\big)} \parn 
The law of the variable $Z_s$ is not at all that simple, but it is known (see [Y], [CDJR]) that its Fourier-Laplace transform is computable. The following lemma is proved in Section \ref{sec.Proofs}. 
\blem  \label{lem.LaplZ} \  The law of the variable $Z_s$ of (\ref{f.2chaoss}) is given by$\,$:  for any $\,s\ge 0$ and real $\,r,c, b$,  
\begin{equation*}
\E_0\bigg[ \exp\!\bigg(\rt1\Big[ r\, w_s + c \int_0^s w_\tau\, d\tau\Big] -  { \frac{b^2}{2}} \int_0^s w_\tau^2\, \, d\tau\bigg)\bigg] 
\end{equation*}
\begin{equation*}
=\,  \frac{1}{\sqrt{\ch\!(b s)}} \, \exp\!\bigg[- \frac{\th\!(b s)}{2\, b} \,r^2 - 2\, \frac{ \sh\!^2(b s/2)}{b^2\, \ch\!(b s)}\, r c - \frac{ {bs} -\th\!(b s)}{2\,b^3}\, c^2 \bigg] . 
\end{equation*}
\elem   \vspace{-2mm} 
Of course, for $\,b=0\,$ we recover the Fourier transform of (\ref{f.Mehl}), namely $\,{e^{-\big(r^2+s\, rc + \frac{s^2}{3} c^2\big) s/2}}$.   
\bpro  \label{pro.LaplZ} \  The $\,x$-Laplace transform of the variable $Z_1$ of (\ref{f.2chaoss}) is given by$\,$:  for any real $\,w, z, b$,  \vspace{-3mm} 
\begin{equation*}
\int_{0}^\infty e^{- {\frac{b^2}{2}} x} \,\, q^0_1(w , z  , x )\, dx\, =\, \frac{b^2\,\exp\!\bigg[  -
\frac{\big[b\, z - \th\!(b/2)\,w\big]^2 +\, \coth b\, \big[b - 2\,\th\!(b /2)\big] w^2}{2\, \big[1 - (2/b)\,\th\!(b /2)\big]} \bigg] }{2\pi \sqrt{\big[b - 2\,\th\!(b /2)\big] \sh b\,}} 
\end{equation*}  \vspace{-2mm} 
\begin{equation*}
= \, \frac{b^2}{2\pi \sqrt{\big[b - 2\,\th\!(b /2)\big] \sh b\,}} \times \exp\!\Bigg[ \frac{b^2}{8}\!\times\!\frac{(w-2z)^2}{1-\frac{b}{2}\,\coth\!(\frac{b}{2})} - \frac{b^2}{2}\, z^2 - {\frac{b}{4}}\,\coth\!({\ts\frac{b}{2}})\,w^2\Bigg] .
\end{equation*}
\epro 
This is of course consistent with (\ref{f.Mehl}), via $\,b\to0\,$; and integrating with respect to $\,dw\,dz\,$, we recover \  ${\ds \int_{0}^\infty e^{- {\frac{b^2}{2}} x} \bigg[\int\!\!\int q^0_1(w , z  , x )\, dw\,dz\bigg] dx = (\ch b)^{-1/2}}$, \  as it must be. \parn 
\ub{Proof} \quad  We invert the Fourier transform in Lemma \ref{lem.LaplZ} by Plancherel's Formula : 
\begin{equation*}
\int_{0}^\infty  e^{- {\frac{b^2}{2}} x} \; q^0_1(w , z  , x )\, dx 
\end{equation*}   \vspace{-3mm} 
\begin{equation*}
=\, \frac{1}{4\pi^2 \sqrt{\ch b}} \int_{\sR^2} e^{-\sqrt{-1}\, [ w\, r + z\, c]} \, \exp\!\bigg[- \frac{\th b}{2\, b} \,r^2 - 2\, \frac{ \sh\!^2(b/2)}{b^2\, \ch b}\, r c  - \frac{ {b} -\th b}{2\,b^3}\, c^2 \bigg] dr\, dc 
\end{equation*}
\begin{equation*}
=\, \frac{1}{4\pi^2 \sqrt{\ch b}} \int_{\sR^2} e^{-\sqrt{-1}\, [ w\, r + z\, c]} \, \exp\!\bigg[- \frac{\th b}{2\, b} \,\Big(r + {\frac{\th\!(b/2)}{b}}\, c\Big)^2 - \frac{ {b} - 2\,\th\!(b /2)}{2\,b^3}\, c^2\bigg] dr\, dc 
\end{equation*}
\begin{equation*}
=\, \frac{1}{4\pi^2 \sqrt{\ch b}} \int_{\sR^2} e^{-\sqrt{-1}\, \big[ w\, r + \big(z - {\frac{\th\!(b/2)}{b}}\,w\big) c\big]} \, \exp\!\bigg[- \frac{\th b}{2\, b} \,r^2 - \frac{ {b} - 2\,\th\!(b /2)}{2\,b^3}\, c^2\bigg] dr\, dc 
\end{equation*}
\begin{equation*}
=\, \sqrt{\frac{b}{\sh b}}\times \frac{e^{- \frac{b\,w^2}{2\, \th b}}}{2\pi}
\int_{\sR} e^{-\sqrt{-1}\, \big(z - {\frac{\th\!(b/2)}{b}}\,w\big) c} \, \exp\!\bigg[ - \frac{ {b} - 2\,\th\!(b /2)}{2\,b^3}\, c^2\bigg] \frac{dc}{\sqrt{2\pi}}
\end{equation*}
\begin{equation*}
=\, \sqrt{\frac{b}{\sh b}}\times \frac{e^{- \frac{b\,w^2}{2\, \th b}}}{2\pi} \times \sqrt{\frac{b^3}{{b} - 2\,\th\!(b /2)}}\times \exp\!\bigg[ - \frac{b^3}{2 \big[b - 2\,\th\!(b /2)\big]}\, \big(z - {\ts\frac{\th\!(b/2)}{b}}\,w\big)^2\bigg]
\end{equation*}
\begin{equation*}
=\, \frac{b^2}{2\pi \sqrt{\big[b - 2\,\th\!(b /2)\big] \sh b\,}} \times \exp\!\Bigg[  -
\frac{\big[b z - \th\!(b/2)\,w\big]^2 +  \coth b\, \big[b - 2\,\th\!(b /2)\big] w^2}{2\, \big[1 - (2/b)\,\th\!(b /2)\big]} \Bigg] . \;\;\diamond 
\end{equation*}

\subsection{Fourier-Laplace transform of the tangent process $\,Y_s$} \label{sec.FLTdensYs} \indf 
   We use Lemma \ref{lem.LaplZ} to express this Fourier-Laplace transform. 
\blem  \label{lem.LaplY} \  The law of the variable $Y_s$ of (\ref{f.2chaos}) is given by$\,$:  for any $\,s\ge 0$ and real $\,r,\rr, b, \gamma , c $,  
\begin{equation*}
\E_0\Big[ \exp\!\Big(\rt1\big[ r\, w_s + \rr\, \beta_s + \gamma \,\zeta_s + c \, \bar z_s \big] -  b^2 A_s\Big)\Big] 
\end{equation*}
\begin{equation*}
=\,  \frac{1}{\ch\!(b s)} \, \exp\!\bigg[- \frac{\th\!(b s)}{2\, b} \,(r^2+\rr^2) - 2\, \frac{ \sh\!^2(b s/2)}{b^2\, \ch\!(b s)}\, (r\, c + \rr\, \gamma) - \frac{ {bs} -\th\!(b s)}{2\,b^3}\, (c^2+\gamma^2) \bigg] . 
\end{equation*}
In particular, the law of $\,A_s\,$ is given by \  $\E_0\big[ \exp( - b^2 A_s)\big] =\, {1}/{\ch\!(b s)}$\raise0.9pt\hbox{.} 
\elem 
\ub{Proof} \quad  This follows directly from Lemma \ref{lem.LaplZ}, by  independence of $\,w\,$ and $\,\beta\,$. 
$\; \diamond $   

\bpro  \label{pro.LaplY} \  The $\,x$-Laplace transform of the variable $Y_1$ of (\ref{f.2chaos}) is given by$\,$:  for any real $\,w,\beta , z , \zeta , b$,  
\begin{equation*}
\int_{0}^\infty  e^{- {{b^2}} x} \; q_1(w,\beta , x , \zeta , z) \, dx
\end{equation*}
\begin{equation}  \label{f.Psi(b)} 
=\, \frac{b^4\, \exp\!\Big[ \frac{b^2}{8}\!\times\!\frac{(w-2z)^2+(\beta-2\zeta)^2}{1-\frac{b}{2}\,\coth\!(\frac{b}{2})} - \frac{b^2}{2}(z^2+\zeta^2) - {\ts\frac{b}{4}}\,\coth\!({\ts\frac{b}{2}})(w^2+\beta^2) \Big]} {8\pi^2\, \big[b\, \ch\!({\ts\frac{b}{2}}) - 2\,\sh\!({\ts\frac{b}{2}})\big]\,  \sh\!({\ts\frac{b}{2}})} \, =: \, \Psi_{w,\beta,\zeta,z}(b)\, \raise0pt\hbox{.}
\end{equation}
In particular, $\,\Psi_{w,\beta,\zeta,z}(0) = \frac{3}{\pi^2} \exp\!\Big[  - \frac{w^2+\beta^2}{2} - 6 \big(z - w/2\big)^2 - 6 \big(\zeta - \beta/2\big)^2\Big]$ is the marginal density of $\big( w_1\, , \, \beta_1\, , \,\zeta_1\, , \,\bar z_1\big)$.  
\epro 
\ub{Proof} \quad  As Lemma \ref{lem.LaplY} follows from Lemma \ref{lem.LaplZ}, this follows merely from Proposition \ref{pro.LaplZ} by  independence of $\,w\,$ and $\,\beta\,$. Indeed, for any test functions $\,f,g\,$ on $\R^2$ we have\,:  
\begin{equation*}
\int_{\sR^4} f(w,z)\, g(\beta,\zeta) \left[\int_{0}^\infty  e^{- {{b^2}} x} \, q_1(w,\beta , x , \zeta , z) \, dx\right] dw\, dz\, d\beta\, d\zeta
\end{equation*}
\begin{equation*}
= \E\Big[ f(w_1,\bar z_1)\, g(\beta_1,\zeta_1)\,  e^{- {\frac{b^2}{2}} \big[\int_0^1w_\tau^2\, d\tau+\int_0^1\beta_\tau^2\, d\tau\big]} \Big] 
\end{equation*}
\begin{equation*}
= \E\Big[ f(w_1,\bar z_1)\,  e^{- {\frac{b^2}{2}} \int_0^1w_\tau^2\, d\tau} \Big] \times \E\Big[ g(\beta_1,\zeta_1) \,  e^{- {\frac{b^2}{2}}\int_0^1\beta_\tau^2\, d\tau} \Big] 
\end{equation*}
\begin{equation*}
= \int_{\sR^2} f(w,z) \int_{0}^\infty  e^{- {\frac{b^2}{2}} x} \, q_1^0(w,z, x) \, dx\,\, dw\, dz\times \int_{\sR^2} g(\beta,\zeta) \int_{0}^\infty  e^{- {\frac{b^2}{2}} x} \, q_1^0(\beta ,\zeta, x) \, dx\,\, d\beta\, d\zeta
\end{equation*}
\begin{equation*}
= \int_{\sR^4} f(w,z)\, g(\beta,\zeta) \left[\int_{0}^\infty  e^{- {\frac{b^2}{2}} x} \, q_1^0(w,z, x) \, dx\times\! \int_{0}^\infty  e^{- {\frac{b^2}{2}} x} \, q_1^0(\beta ,\zeta, x) \, dx\right] dw\, dz\, d\beta\, d\zeta\, , 
\end{equation*}
and the claim follows directly from Proposition \ref{pro.LaplZ}. 
\if{   
We invert the Fourier transform in Lemma \ref{lem.LaplY}, as for Proposition \ref{pro.LaplZ} : 
\begin{equation*}
\int_{0}^\infty  e^{- {{b^2}} x} \;  q_1(w,\beta , x , \zeta , z)\, dx\,
=\, \frac{1}{16\pi^4 \,\ch b} \int_{\sR^4} e^{-\sqrt{-1}\, [ w\, r + z\, c + \beta\, t + \zeta\, u]}\, \times 
\end{equation*}
\begin{equation*}
\exp\!\Bigg[- \frac{\th b}{2\, b} \bigg[\Big(r + {\frac{\th\!(b/2)}{b}}\, c\Big)^2 + \Big(t + {\frac{\th\!(b/2)}{b}}\, u\Big)^2\bigg] - \frac{ {b} - 2\,\th\!(b /2)}{2\,b^3}\, (c^2+u^2)\Bigg] dr\, dc \, dt\, du
\end{equation*}
\begin{equation*}
=\, \frac{1}{16\pi^4 \,\ch b} \int_{\sR^4} e^{-\sqrt{-1}\, \big[ w\, r + \big(z - {\frac{\th\!(b/2)}{b}}\,w\big) c\big] -\sqrt{-1}\,  \big[ \beta\, t + \big(\zeta - {\frac{\th\!(b/2)}{b}}\,\beta\big) u\big]}\, \times \qquad \qquad $$
$$ \qquad \qquad \qquad \qquad \times \exp\!\bigg[- \frac{\th b}{2\, b} (r^2+t^2) - \frac{ {b} - 2\,\th\!(b /2)}{2\,b^3}\, (c^2+u^2)\bigg] dr\, dc \, dt\, du\, , 
\end{equation*}
which splits into two terms of the kind already seen for Proposition \ref{pro.LaplZ}. This yields 
\begin{equation*}
\int_{0}^\infty  e^{- {{b^2}} x} \; q_1(w,\beta , x , \zeta , z) \, dx\; = \; \frac{b^4}{4\pi^2 \big[b - 2\,\th\!(b /2)\big]\, \sh b\,}\, \times \hskip 30mm   $$
$$ \hskip 11mm  \times\, \exp\!\Bigg[  - \, 
\frac{\big[b\, z - \th\!(b/2)\,w\big]^2 + \big[b\, \zeta - \th\!(b/2)\,\beta\big]^2 + \coth b\, \big[b - 2\,\th\!(b /2)\big] (w^2+\beta^2)}{(2/b) \big[b - 2\,\th\!(b /2)\big]} \Bigg] 
\end{equation*}
\begin{equation*}
= \frac{b^4}{4\pi^2 \big[b - 2\,\th\!(b /2)\big] \sh b} \exp\!\Bigg[  
\frac{[b\,\coth b-1] (w^2+\beta^2) - 2b\,\th\!(b/2)\,[wz+\beta\zeta] + b^2[z^2+\zeta^2]}{-\, 2\, \big[1 - (2/b)\,\th\!(b /2)\big]} \Bigg] 
\end{equation*}
\begin{equation*}  
= \frac{b^4\, e^{\frac{b^2}{8} \big[w^2+\beta^2-4wz-4\beta\zeta\big]}} {8\pi^2\big[b\, \ch\!({\ts\frac{b}{2}}) - 2\,\sh\!({\ts\frac{b}{2}})\big]  \sh\!({\ts\frac{b}{2}})} \, \exp\!\Bigg[ \frac{b^2\big[(w-2z)^2+(\beta-2\zeta)^2\big]}{8\big[\frac{2}{b}\,\th\!(\frac{b}{2}) -1\big]} -  {\ts\frac{b}{4}}\,\coth\!({\ts\frac{b}{2}})(w^2+\beta^2) \Bigg] , 
\end{equation*}
and the result follows directly. 
}\fi 
$\;\diamond$ \par 

\blem  \label{lem.zthz} \  All solutions $\,z\in\C^*$ of the equation $\,z=\th z\,$ belong to the imaginary axis, and form a sequence $\,\RR= \{\pm\rt1 y_n\,|\,n\in\N\}$, with $\,\frac{17\pi}{12} < y_0 <y_1<\ldots < y_n\nea\infty\,$. 
\elem 
\ub{Proof} \  Writing $\,z=x+\rt1 y \in\C^*$, we have $\,z=\th z\, \LRa\, e^{2z} = \frac{1+z}{1-z}\,$ and then equivalently 
\begin{equation*}
\cos(2y) =\,  {\ts\frac{1-x^2-y^2}{(1-x)^2+y^2}}\, e^{-2x} \quad \hbox{ and }  \quad \sin(2y) =\,  {\ts\frac{2\,y}{(1-x)^2+y^2}}\, e^{-2x}\,\, ,
\end{equation*}
whence 
\begin{equation*}
e^{4x}\big((1-x)^2+y^2\big)^2 = 4y^2 + (1-x^2-y^2)^2\, . 
\end{equation*}
The latter is equivalent either to $\,x=0\,$, or to $\,y^2= -(1-x)^2$ (which is excluded), or to \  $y^2 = \frac{(1+x)^2-e^{4x} (1-x)^2}{e^{4x}-1}\, $\raise1.7pt\hbox{.}   \    Then using this last value of $\,y^2$, by the above we must also have \  ${ \sin(2y) = \frac{y}{x}\, \sh(2x)}$, 
which is impossible since for any $\,x,y\in\R^*$ we have \  ${ \frac{\sin(2y)}{y} <2<  \frac{\sh(2x)}{x}\,}$\raise1.3pt\hbox{,}  and clearly $\,z\,$ cannot be real. 
Hence we are left with $\,z = \pm\rt1 y\,$, with \  $y>0\,$ and then $\,y=\tg y\,$. The claim follows, with moreover $\big((n+ {3}/{2})\pi - y_n\big) \sea 0\,$, and $\, \frac{3\pi}{2} > y_0> \frac{17\pi}{12}\,$ since $\,\tg \frac{17\pi}{12} = \cotg \frac{\pi}{12} =2+\sqrt{3} < \frac{17\pi}{12}\,$\raise1.7pt\hbox{.} $\;\diamond $  \parm

   The Laplace transform $\,\Psi_{w,\beta,\zeta,z}(b)$ in Proposition \ref{pro.LaplY}
is a meromorphic function of $\,b\in\C\,$, with singularities  at the points of $\,\rt1 2\pi\,\Z^*$ and at the non-null zeros of  $\big[b - 2\,\th\!(b /2)\big]$, that is to say at the points of $\,2\RR$ (according to Lemma \ref{lem.zthz}). 
\par 

\if{ 
\ub{Proof} \quad Let $\sqrt{b}\,>0$ on $\R_+^*\,$.  The function  $\,b\mapsto \Psi_{w,\beta,\zeta,z}\big(\sqrt{b}\,\big)\,$ is analytic in $\C\moins ]-\infty,0]\,$.  
According to Proposition \ref{pro.LaplY} we have 
\begin{equation*}
\int_{0}^\infty  e^{- {{b}}\, x} \;  q_1(w,\beta , x , \zeta , z)\, dx\,
=\, \Psi_{w,\beta,\zeta,z}\big(\sqrt{b}\,\big)\, , 
\end{equation*}
and then for any $\,\eta > 0\,$: \vspace{-2mm} 
\begin{equation*}
q_1(w,\beta , x , \zeta , z) = \frac{-\rt1}{2\pi} \int_{\eta -\rt1\infty}^{\eta +\rt1\infty}  e^{x\, b} \; \Psi_{w,\beta,\zeta,z}\big(\sqrt{b}\,\big)\, db 
\end{equation*}
\begin{equation}  \label{f.exprp1}
 = \frac{e^{\eta\,x}}{2\pi} \int_{-\infty}^{\infty}  e^{\rt1 x\,y}\,\, \Psi_{w,\beta,\zeta,z}\Big(\sqrt{\eta +\rt1 y}\,\Big)\, dy\, . 
\end{equation}
Then the continuity of $\, \Psi_{w,\beta,\zeta,z}\big(\sqrt{b}\,\big)\,$ at $\,b=0\,$ allows to let $\,\eta\sea 0\,$, yielding 
\begin{equation*}
q_1(w,\beta , x , \zeta , z) = 
\end{equation*}
\begin{equation*}
\frac{1}{2\pi} \Bigg[\!\int_{0}^{\infty}\! e^{\rt1 x\,y}\, \Psi_{w,\beta,\zeta,z}\big((1+ \rt1\!)\sqrt{y/2}\,\big) dy + \int_{-\infty}^{0}\! e^{\rt1 x\,y}\, \Psi_{w,\beta,\zeta,z}\big((1- \rt1\!)\sqrt{-y/2}\,\big) dy \Bigg]
\end{equation*}
\begin{equation*}
= \frac{1}{\pi} \int_{0}^{\infty}  \Re\Big[e^{\rt1 x\,y}\, \Psi_{w,\beta,\zeta,z}\big((1+ \rt1\!)\sqrt{y/2}\,\big) \Big] dy \, , 
\end{equation*}
from which Formula (\ref{f.p1=}) follows, provided $\,y\mapsto \Big| \Psi_{w,\beta,\zeta,z}\big((1+ \rt1\!) \sqrt{y}\,\big)\Big| $ be integrable on $\R_+\,$. \   Now we have \  $ \sh\!\big((1+ \rt1\!)\, y\big) = \sh y\,\cos y + \rt1 \ch y\,\sin y \; ;  $ 
$$ \th\!\big((1+ \rt1\!)\, y/2\big) = {\ts\frac{\sh y + \rt1 \sin y}{\ch y + \cos y}}  = 1 +\O(e^{-y})\; ;\; \coth\!\big((1+ \rt1\!)\, y/2\big) = {\ts\frac{\sh y - \rt1 \sin y}{\ch y - \cos y}} \,\raise1.5pt \hbox{.}  $$ 
Hence by Proposition \ref{pro.LaplY}, as $\,y\to\infty\,$ we have 
$$ \Big| \Psi_{w,\beta,\zeta,z}\big((1+ \rt1\!)\, y\big)\Big| =  \O(y^3\, e^{-y}) \exp\!\Big[  
- {\ts \frac{w^2+ \beta^2}{4}}\, \big[y +\O(1)\big] \Big] , $$
which shows the wanted integrability, and then Formula (\ref{f.p1=}). $\;\diamond$ \par 
}\fi  

{\small   
\subsection{Complement\,: Density $\,\alpha_s(x)$ of the variable $\,A_s\,$} \indf 
   Denote by $\,\alpha_s = \alpha_s(x)\,$ the density of the variable $\,A_s\,$, so that for any $\,s>0\,$ we have \parn
\centerline{${\ds\, \alpha_s(x) = \int q_s(w,\beta , x , \zeta , z)\, dw\,d\beta\, d\zeta\, dz\,=\, \frac{1}{s^2}\, \alpha_1\left(\frac{x}{s^2}\right)}$, \   by (\ref{f.densYs}).} 
\blem  \label{lem.propralpha} \   The density $\,\alpha_1\,$ is smooth and bounded (with bounded derivatives), and we have \   
\begin{equation*}  
\alpha_1(x) =\, \frac{4}{\pi} \int_{0}^{\infty}\,  \frac{\cos\!\big(2x\,y^2 - y\big) } { \sh\!^2 y + \cos^2 y } \; y\, \sh y\, dy + \frac{4}{\pi} \int_{0}^{\infty}\,  \frac{\cos y\, \cos(2x\,y^2)} { \sh\!^2 y + \cos^2 y } \; e^{-y}\, y\, dy \,  .
\end{equation*}
\elem 
\ub{Proof}  \quad   According to Lemma \ref{lem.LaplY}, 
for any $\,\eta > 0\,$ we have\,: \quad $\alpha_1(x) =\, $ 
\begin{equation*}   
\frac{e^{\eta\,x}}{2\pi} \int_{-\infty}^{\infty}  \frac{e^{\rt1 x\,y}}{\ch\! \Big(\sqrt{\eta +\rt1 y}\,\Big)}\,  dy\,  
=\, \frac{e^{\eta\,x}}{2\pi} \int_{-\infty}^{\infty}  \frac{e^{\rt1 x\,y}}{\ch\! \bigg(\sqrt{\frac{\sqrt{y^2+\eta^2} +\eta}{2}} + \rt1 {\rm sgn}(y) \sqrt{\frac{\sqrt{y^2+\eta^2} -\eta}{2}}\,\,\bigg)}\;  dy\,   .
\end{equation*}
Expanding and letting  $\,\eta \sea 0\,$ we obtain$\,$: 
\begin{equation*}   
\alpha_1(x) =\, \frac{4}{\pi} \int_{0}^{\infty}\,  \frac{ \ch y\, \cos y\, \cos(2x\,y^2) + \sh y\, \sin y\, \sin(2x\,y^2) } { \sh\!^2 y + \cos^2 y } \; y\, dy \,  
\end{equation*}
\begin{equation*}  
=\, \frac{4}{\pi} \int_{0}^{\infty}\,  \frac{\cos\!\big(2x\,y^2 - y\big) } { \sh\!^2 y + \cos^2 y } \; y\, \sh y\, dy + \frac{4}{\pi} \int_{0}^{\infty}\,  \frac{\cos y\, \cos(2x\,y^2)} { \sh\!^2 y + \cos^2 y } \; e^{-y}\, y\, dy \,  .
\end{equation*}
As a consequence,  $\,\alpha_1\,$ is smooth, and bounded by $\,{\ds \frac{4}{\pi} \int_{0}^{\infty} \frac{y\, \ch y\; dy} { \sh\!^2 y + \cos^2 y } }\,$\raise1.8pt\hbox{.}   $\;\diamond$ 
   
\brem  \label{lem.denonul} \   We have \   $\,\sh\!^2 z + \cos^2\!z = 0\,\LRa\, \ch\!^2 z = \sin^2\!z \,\LRa\, z\in (1\pm\rt1\!) {\ds \frac{\pi}{4}} ( 1+ 2\Z)$,  \parn 
\centerline{ i.e.,  \  $ |z|\in \frac{\pi ( 1+ 2\sN)}{2\sqrt{2}}\,\raise1.pt\hbox{,} \quad {\rm Arg}\, z\, \equiv\, \frac{\pi}{4}\,$ modulo $\,\frac{\pi}{2}\,$\raise1.8pt\hbox{.} }
\erem 
\ub{Proof}  \quad There is clearly no real solution, since $\,\ch x> 1\ge \sin x\,$ for any $\,x\in\R^*$. Consider then $\,z=(\rt1 + \alpha)\, x\,$, with $\,x,\alpha\in \R\,$. Up to change $\,z\,$ into $-z\,$, we only have to consider the equation $\,\ch z = \sin z\,$. Then 
$$ \ch z = \sin z \;\LRA\; \ch\!(\alpha x) \cos x = \ch x\, \sin(\alpha x) \; \hbox{ and }\; \sh\!(\alpha x) \sin x = \sh x\, \cos(\alpha x) $$ 
$$ \Ra\, \ch\!^2(\alpha x) \cos^2x +\sh\!^2(\alpha x) \sin^2\!x = \ch\!^2x\, \sh^2x \LRA \big(\ch\!^2x - \cos^2\!x\big)\big(\sh\!^2(\alpha x) - \sh\!^2x\big) = 0\, .$$ 
As $\,x=0\,$ is not a solution, the only possibility is \  $\sh\!(\alpha x) = \pm\, \sh x\,$, whence \  $\alpha = \pm1$, and then $\,\cos x= \alpha\, \sin x\,$.  This yields the claim. $\;\diamond$ 
\bpro \label{pro.densY} \  We have \quad    ${\ds \alpha_1(x) =\, 2\pi\sum_{n\in\sN}\, (-1)^n \big(n+\5\big)\, e^{-\big(n+\5\big)^2\pi^2\, x} }\,$, \  so that  \parn   
\centerline{${\ds \pi\, e^{-\pi^2\, x/4} \big(1- 3\, e^{-2\pi^2\, x} \big) \le \, \alpha_1(x) \le\, \pi\, e^{-\pi^2\, x/4} }$ \   for any $\,x\ge \pi\2$, \  and then} \parn
\centerline{${\ds \alpha_1(x)\,= \, \pi\, e^{-\pi^2\, x/4}\, \big(1- \O( e^{-2\pi^2\, x}) \big) }\,$ as $\,x\to\infty$\,. Moreover, $\,\alpha_1\,$ decreases on $[3\pi\2,\infty[\,$.} \parn 
\epro
\ub{Proof} \quad  This results from Lemma \ref{lem.LaplY} and \big([BPY], Table 1 continued and Section 3.3, in particular Formula (3.11), with $\,C_1\equiv 2A_1$, got merely by expanding $\,{1}/({\ch t})$\big). \   
For any $\,\pi^2x> 1$, this alternate series has decreasing generic term  ${\ds \big(n+\5\big)\, e^{-\big(n+\5\big)^2\pi^2\, x} }\,$, whence the estimate. Finally, the same holds for the series ${\ds \big(n+\5\big)^3\, e^{-\big(n+\5\big)^2\pi^2\, x} }\,$ yielding $\,\alpha_1'(x)$, thereby garanteing $\,\alpha_1'(x)<0$, as soon as $\,\pi^2x> 3$. $\;\diamond$ \parsn

 Lemma \ref{lem.LaplY} and Proposition \ref{pro.densY} at once entail the following. \vspace{-2mm} 
\bcor \label{cor.densY} \  ${\ds \lambda \lmt \int_0^\infty e^{\lambda\, x}\, \alpha_1(x)\, dx}\,$ defines an analytic function on $\,\Re(\lambda) < \pi^2/4\,$, equal to    $\,{ \frac{1}{\cos\sqrt{\lambda}} }\,$ for $\,0\le \lambda< \frac{\pi^2}{4}\,$\raise1.5pt\hbox{.} Moreover, we have \  $\alpha_s(x) = \frac{\pi}{s^2}\,e^{-\pi^2x/(4s^2)} \big(1- \O( e^{-2\pi^2\, x/s^2}) \big)$ (as $\,\frac{x}{s^2} \to\infty$). 
\ecor 
\ub{Proof} \quad  By  Lemma \ref{lem.LaplY} the set of $\,b\in\C$ such that \  ${\ds \int_0^\infty e^{-b^2\, x}\, \alpha_1(x)\, dx = \frac{1}{\ch b} }\,$ contains $\R\,$. By Proposition \ref{pro.densY}  \  ${\ds b \mapsto \int_0^\infty e^{-b^2\, x}\, \alpha_1(x)\, dx}\,$ is analytic on $\big\{|\Im(b)|<\pi/2\big\}$, and so is ${\ds b \mapsto  \frac{1}{\ch b} }\,$ too. Hence we have \  ${\ds \int_0^\infty e^{b^2\, x}\, \alpha_1(x)\, dx = \frac{1}{\cos b} }\,$ for $\,|\Re(b)|<\pi/2\,$.  $\;\diamond$

\brem \label{rem.loiA1} \  {\rm  1) \  The law of $\,2\,A_1\,$ is that of \  $\,\inf\big\{\tau>0\,\big|\, |\beta_\tau| = 1\big\}$ \   and also that of \parn 
$\Big(\max\big\{|\beta_\tau|\,\big|\, 0\le \tau\le 1\big\}\Big)\2$ \big(See [BPY], Table 2, which also exhibits two other random variables having the same law as $\,A_1$\big).  \quad 
2) \  Note that for $\,\Re(\lambda) < \pi^2/4\,$ (letting $\sqrt{\cdot}\,>0\,$ on $\R_+^*$) : 
\begin{equation*}
\frac{1}{\cos\sqrt{\lambda}} = \int_0^\infty\!\! \int  e^{\lambda\, x}\, q_1(w,\beta , x , \zeta , z)\, dw\,d\beta\, d\zeta\, dz\, dx\, =  \int  \Psi_{w,\beta,\zeta,z}\big(\sqrt{-\lambda}\,\big) \, dw\,d\beta\, d\zeta\, dz\, .
\end{equation*}
}\erem 
}  
\vspace{-5mm}

\section{Integral expression of the density $\,q_s(w,\beta , x , \zeta , z)$} \label{s.intexprdqs}  
   Let \vspace{-3mm} 
\begin{equation} \label{f.qtildes} 
\tilde q_s\equiv \tilde q_s(w,\beta , \zeta , z) :=\, \frac{3}{\pi^2 s^4}\, \exp\!\bigg[- { \frac{3  \big[(sw-2z)^2+(s\beta-2\zeta)^2\big] }{2\,s^3}} - {\frac{w^2+\beta^2}{2\, s}} \bigg]  
\end{equation}
denote the marginal Gaussian density of  \  ${\ds \tilde Y_{s}  \,:= \,\Big( w_s\, , \;  \beta_s\, , \int_0^s \beta_\tau\,d\tau\, ,  \int_0^s w_\tau\,d\tau \Big)}$ (according to (\ref{f.Mehl}) or Proposition \ref{pro.LaplY})). \parsn 
\ub{\bf Notation} \quad   Consider the function  \  $\,\Phi(\lambda)\equiv \Phi_{w,\beta,\zeta,z}(\lambda) :=\, \Psi_{w,\beta,\zeta,z}\big(\sqrt{-\lambda}\,\big)$, derived from Proposition \ref{pro.LaplY}. We systematically use the usual determination of the complex square root, cutting $\C$ along the negative real semi-axis and letting $\sqrt{\cdot}\,>0\,$ on $\R_+^*$. \   By the expression (\ref{f.Psi(b)}), we have\,:  
\begin{equation} \label{f.Psi-lambd} 
\Phi_{w,\beta,\zeta,z}(\lambda) =\, \frac{\lambda^2\, \exp\!\Big[\frac{\lambda}{8}\!\times\! \frac{(w - 2z)^2 + (\beta - 2\zeta)^2}{\frac{\sqrt{\lambda}}{2}\,\cotg\!(\frac{\sqrt{\lambda}}{2}) - 1} + \frac{\lambda}{2}(z^2+\zeta^2) - (w^2+\beta^2)\frac{\sqrt{\lambda}}{4}\,\cotg\!({\frac{\sqrt{\lambda}}{2}}) \Big]}{8\pi^2\, \big[2 \sin(\frac{\sqrt{\lambda}}{2}) - \sqrt{\lambda} \cos(\frac{\sqrt{\lambda}}{2})\big] \sin(\frac{\sqrt{\lambda}}{2})} \, \raise1.9pt\hbox{.}
\end{equation}

\bLem  \label{lem.explPsi} \  The function \  $\lambda\mapsto \Phi_{w,\beta,\zeta,z}(\lambda)\,$ is analytic for $\,\Re(\lambda)< 4\pi^2$.
\eLem \vspace{-2mm} 
\ub{Proof} \quad  Note that the functions  \   
$\sin(\frac{\sqrt{\lambda}}{2})\Big/\frac{\sqrt{\lambda}}{2}\,$ and $\,\cos(\frac{\sqrt{\lambda}}{2})\,$ are plainly analytically continued for any $\,\lambda\in\C$, and that by the expression (\ref{f.Psi-lambd}), $\Phi(\lambda)$ is analytically continued at $\,\lambda =0$, and analytic at any $\,\lambda\in\C^*$ such that $\,\sin(\frac{\sqrt{\lambda}}{2}) \not = 0\,$ and $\,\tg\!(\frac{\sqrt{\lambda}}{2}) \not = \frac{\sqrt{\lambda}}{2}\,$\raise1.1pt\hbox{,} hence,  according to Lemma \ref{lem.zthz}, at those $\lambda\in\C$ not belonging to the sequence $\{ 4\pi^2 (n+1)^2, 4\, y_n^2\, \,|\,n\in\N\} \subset\, [4\pi^2,\infty[\,$. \   
This shows the analyticity for $\,\Re(\lambda) < 4\pi^2$.   $\;\diamond$ \parsn   
The following is proved in Section \ref{sec.Proofs}.   
\bPro \label{pro.abscCv} \quad   For any $\,(w,\beta,\zeta,z)\,$ we have
\begin{equation*}
\Phi_{w,\beta,\zeta,z}(\lambda) = \int_{0}^\infty  e^{\lambda\, x} \; q_1(w,\beta , x , \zeta , z) \, dx\;,  \quad \hbox{ for } \;  \Re(\lambda) < 4\pi^2, \hbox{ and} 
\end{equation*}
\begin{equation}   \label{f.invFourq1} 
e^{r\, x} \, q_1(w,\beta , x , \zeta , z) = \int_{-\infty}^\infty  e^{-\rt1 t\, x}\, \Phi(r +\rt1 t) \, \frac{dt}{2\pi}\, \raise1.7pt\hbox{,} \quad \hbox{ for } \, x>0  \, \hbox{ and } \; r< 4\pi^2 . 
\end{equation}
\ePro \vspace{-2mm} 
By scaling and using (\ref{f.invFourq1}), for any $\,s\, ,\, x\, >0\,$, $r<4\pi^2$  and  $(w,\beta , \zeta , z)\in\R^4\,$ we have\,:  
\begin{equation*}  
q_s(w,\beta , x , \zeta , z)\, =\, \frac{1}{s^6}\,\, q_1\bigg( \frac{w}{\sqrt{s}}\, \raise1.8pt\hbox{,}\, \frac{\beta}{\sqrt{s}}\, \raise1.8pt\hbox{,}\,  \frac{x}{s^2}\, \raise1.8pt\hbox{,}\, \frac{\zeta}{\sqrt{s^3}}\, \raise1.8pt\hbox{,}\, \frac{z}{\sqrt{s^3}} \bigg) 
\end{equation*}
\begin{equation*}  
=\, \frac{e^{- r\, x/s^2}}{2\pi\, s^6} \int_{-\infty}^\infty  e^{-\rt1 x\,y/s^2}\; \Phi_{\big( \frac{w}{\sqrt{s}}\, \raise1.8pt\hbox{,}\, \frac{\beta}{\sqrt{s}}\, \raise1.8pt\hbox{,}\,  \frac{\zeta}{\sqrt{s^3}}\, \raise1.8pt\hbox{,}\, \frac{z}{\sqrt{s^3}} \big)}(r + \rt1 y)\, dy\,  . 
\end{equation*}
Taking merely  $\, r =0$, \  and setting \quad   $\widetilde \Psi_s := e^{-\rt1 x\,y/s^2}\, \Phi_{\big( \frac{w}{\sqrt{s}} \raise1.8pt\hbox{,} \frac{\beta}{\sqrt{s}} \raise1.8pt\hbox{,}  \frac{\zeta}{\sqrt{s^3}} \raise1.8pt\hbox{,} \frac{z}{\sqrt{s^3}} \big)}\! (\rt1 y)\,$ \parn
for convenience, \   for any $\,\,s , x\, >0\,$ and $(w,\beta , \zeta , z)\in\R^4\,$ we have\,:  
\begin{equation}  \label{f.exprq'}
q_s(w,\beta , x , \zeta , z)\, =\, \frac{1}{2\pi\, s^6} \int_{-\infty}^\infty \widetilde \Psi_s(y) \, dy\,  =\, \frac{1}{2\pi\, s^6} \int_{0}^\infty \big(\widetilde \Psi_s(y) + \widetilde \Psi_s(-y)\big) \, dy\, .    
\end{equation}
Now according to (\ref{f.Psi-lambd}) we have\,: \quad $\Phi_{\big( \frac{w}{\sqrt{s}} \raise1.8pt\hbox{,} \frac{\beta}{\sqrt{s}} \raise1.8pt\hbox{,}  \frac{\zeta}{\sqrt{s^3}} \raise1.8pt\hbox{,} \frac{z}{\sqrt{s^3}} \big)}\! (\rt1 y)\, = $ 
\begin{equation*}  
\frac{- y^2\, e^{B'_s\, \rt1 y\,s\2 - \frac{\sqrt{\rt1 y}}{4s}\,\cotg\!\big({\frac{\sqrt{\rt1 y}}{2}}\big)(w^2+\beta^2)}}{8\pi^2 \big[1- \cos({\sqrt{\rt1 y}}\,) - (\sqrt{\rt1 y}/2) \sin({\sqrt{\rt1 y}}\,)\big]} \, \exp\!\left[\frac{B^2_s\, \rt1 y\, s^{-3}}{1- \frac{2}{\sqrt{\rt1 y}}\,\tg\!(\frac{\sqrt{\rt1 y}}{2})}\right] ,
\end{equation*}
in which we have set  
\begin{equation}  \label{f.defB2B'}
B^2_s := \frac{(sw - 2z)^2  +  (s\beta - 2\zeta)^2}{8} \quad \hbox{ and  } \quad B'_s :=   \frac{4wz + 4\beta\zeta - s(w^2+\beta^2)}{8}\,\raise1.7pt\hbox{.} 
\end{equation}
\vspace{-4mm} 
\if{\small 
\bRem  \label{rem.developB2B'}  \  {\rm   We have \  $B'_s = \frac{z^2+\zeta^2 - 2 B^2_s}{2\, s}$\raise1.1pt\hbox{,} \  and setting \  $\rr \equiv \rr(w,\beta, \zeta,z) :=\frac{w^2+\beta^2}{z^2+\zeta^2}$\raise1.3pt\hbox{,}  \    we also have\,:  
\begin{equation*}  
B_s^2= {\ts\frac{z^2+\zeta^2}{2}}\big[1+\O(s\sqrt{\rr}+s^2 \rr )\big]\, \  \hbox{ and } \  B'_s = (z^2+\zeta^2)\, \O( \sqrt{\rr}+s \rr )\, . 
\end{equation*}
}\eRem 
}\fi  

   Let us now write out a more tractable expression of $\,\widetilde \Psi_s\,$ introduced above. \parsn
First,  for any real $\,y, \t\,$, setting $\,\xi := \sqrt{\frac{|y|}{2}}\,$\raise1.3pt\hbox{,} \   we successively have\,:  
{\small 
\begin{equation*}
{\ts\sqrt{\rt1 y}} = (1+{\rm sgn}(y)\rt1\!)\,  \xi\; ; \quad \frac{2}{\sqrt{\rt1 y}} = (1-{\rm sgn}(y)\rt1\!)\big/\xi\; ; 
\end{equation*} 
\begin{equation*}
\tg\!\big((1+\rt1\!)\,\t\big) = \frac{\sin(2\t) + \rt1 \sh\!(2\t)}{\ch\!(2\t)+ \cos(2\t)}\; ; \; \cotg\!\big((1+\rt1\!)\,\t\big) = \frac{\sin(2\t) - \rt1 \sh\!(2\t)}{\ch\!(2\t) - \cos(2\t)} \; ; 
\end{equation*}
\begin{equation*}  {\ts
\tg\!\Big(\frac{\sqrt{\rt1 y}}{2}\,\Big) } = \frac{\sin\xi + {\rm sgn}(y) \rt1 \sh\xi}{\ch\xi + \cos\xi}\; ; \
{\ts \cotg\!\Big(\frac{\sqrt{\rt1 y}}{2}\,\Big) } = \frac{\sin\xi - {\rm sgn}(y) \rt1 \sh\xi}{\ch\xi - \cos\xi}\; ;  
\end{equation*}
\begin{equation*}  {\ts 
\frac{2}{\sqrt{\rt1 y}} \, \tg\!\Big(\frac{\sqrt{\rt1 y}}{2}\,\Big) } =\, \frac{(\sh\xi+\sin\xi) + {\rm sgn}(y) \rt1 (\sh\xi-\sin\xi) }{ \xi\,(\ch\xi + \cos\xi)}\; ; 
\end{equation*}
\begin{equation*}
\frac{1}{1- \frac{2}{\sqrt{\rt1 y}} \, \tg\!\Big(\frac{\sqrt{\rt1 y}}{2}\,\Big)} =\, \frac{ \xi\,(\ch\xi + \cos\xi) - (\sh\xi+\sin\xi) + {\rm sgn}(y) \rt1\! (\sh\xi-\sin\xi) }{\xi\,(\ch\xi + \cos\xi) - 2 (\sh\xi+\sin\xi) + 2\xi\1 (\ch\xi - \cos\xi)}\; ;  
\end{equation*}
\begin{equation*}
\frac{\rt1 y}{1- \frac{2}{\sqrt{\rt1 y}} \, \tg\!\Big(\frac{\sqrt{\rt1 y}}{2}\,\Big)} =\, \frac{y\,\xi \left[\rt1\! \left(\xi - \frac{\sh\xi+\sin\xi}{\ch\xi + \cos\xi}\right) - {\rm sgn}(y)\! \left(\frac{\sh\xi-\sin\xi}{\ch\xi + \cos\xi}\right)\right] }{ \left(\xi - \frac{\sh\xi+\sin\xi}{\ch\xi + \cos\xi}\right)^2 + \left(\frac{\sh\xi-\sin\xi}{\ch\xi + \cos\xi}\right)^2} =: U(y)\, ;  
\end{equation*}
\begin{equation*}
\frac{\sqrt{\rt1 y}}{2}\, \cotg\!\Big(\frac{\sqrt{\rt1 y}}{2}\,\Big) = \frac{(\sh\xi+\sin\xi) - {\rm sgn}(y) \rt1 (\sh\xi-\sin\xi) }{ 2\,\xi\1\,(\ch\xi - \cos\xi)} =: V(y)\, ; 
\end{equation*}
\begin{equation*}
\cos\!\Big({\ts\sqrt{\rt1 y}}\,\Big) = \ch\xi\,\cos\xi - {\rm sgn}(y) \rt1 \sh\xi\,\sin\xi \; ; 
\end{equation*}
\begin{equation*}
\sin\!\Big({\ts\sqrt{\rt1 y}}\,\Big) = \ch\xi\,\sin\xi + {\rm sgn}(y) \rt1 \sh\xi\,\cos\xi\; ; 
\end{equation*}
\begin{equation*}
1- \cos\!\Big({\ts\sqrt{\rt1 y}}\,\Big) - \Big({\ts\sqrt{\rt1 y}}/2\Big) \sin\!\Big({\ts\sqrt{\rt1 y}}\,\Big) =\, 
\end{equation*}
\begin{equation*}
1- \ch\xi \cos\xi - {\ts\frac{\xi}{2}} (\ch\xi \sin\xi -\sh\xi \cos\xi) + {\rm sgn}(y) \rt1\! \left[\sh\xi \sin\xi  - {\ts\frac{\xi}{2}}(\ch\xi \sin\xi +\sh\xi \cos\xi)\right]\! .
\end{equation*}
} 
Set 
\begin{equation*}
F(y) := \frac{-1}{1- \cos\big({\sqrt{\rt1 y}}\,\big) - \big(\sqrt{\rt1 y}/2\big) \sin\big({\sqrt{\rt1 y}}\,\big)} \, =\, F_r(\xi) + {\rm sgn}(y)\rt1 F_i(\xi)\, , 
\end{equation*}
with 
\begin{equation*}
F_r(\xi) := \frac{\ch\xi \cos\xi + {\frac{\xi}{2}} (\ch\xi \sin\xi -\sh\xi \cos\xi) -1 }{ (\ch\xi-\cos\xi)\big[(\ch\xi-\cos\xi) - \xi\, (\sh\xi+\sin\xi) +\frac{\xi^2}{2}(\ch\xi+\cos\xi)\big]}\,  \raise1.9pt\hbox{,}
\end{equation*}
and 
\begin{equation*}
F_i(\xi) := \frac{\sh\xi \sin\xi  - { \frac{\xi}{2}}(\ch\xi \sin\xi +\sh\xi \cos\xi) }{ (\ch\xi-\cos\xi)\big[(\ch\xi-\cos\xi) - \xi\, (\sh\xi+\sin\xi) +\frac{\xi^2}{2}(\ch\xi+\cos\xi)\big]}\, \raise1.9pt\hbox{.}
\end{equation*}
Hence 
\begin{equation}  \label{f.expresPsi} 
\widetilde \Psi_s \equiv \widetilde \Psi_s(y) \, =\, \frac{y^2\,F(y)}{8\pi^2}\, \exp\!\left[ \rt1 \frac{B'_s-x}{s^2}\, y -\frac{w^2+\beta^2}{2s}\, V(y) + \frac{B^2_s}{s^3}\,U(y)\right] 
\end{equation}
\begin{equation*}
=\, \frac{y^2\,F(y)}{8\pi^2}\, \exp\!\left[ \rt1\!\! \left(\frac{B^2_s}{s^3}\,U_i(\xi) + \frac{B'_s-x}{s^2} + \frac{\ss w^2+\beta^2}{2s}\, V_i(\xi)\!\right) y\, - \frac{B^2_s}{s^3}\,U_r(\xi) -\frac{\ss w^2+\beta^2}{2s}\, V_r(\xi) \right] ,
\end{equation*}
with
\begin{equation*}
U_i(\xi) :=\, \frac{\xi  \left(\xi - \frac{\sh\xi+\sin\xi}{\ch\xi + \cos\xi}\right) }{ \left(\xi - \frac{\sh\xi+\sin\xi}{\ch\xi + \cos\xi}\right)^2 + \left(\frac{\sh\xi-\sin\xi}{\ch\xi + \cos\xi}\right)^2} \; ;  
\quad 
U_r(\xi) :=\, \frac{2\,\xi^3 \left(\frac{\sh\xi-\sin\xi}{\ch\xi + \cos\xi}\right) }{ \left(\xi - \frac{\sh\xi+\sin\xi}{\ch\xi + \cos\xi}\right)^2 + \left(\frac{\sh\xi-\sin\xi}{\ch\xi + \cos\xi}\right)^2} \; ;  
\end{equation*}
\begin{equation*}
V_i(\xi):= \frac{ \sh\xi-\sin\xi }{ 4\,\xi\, (\ch\xi - \cos\xi)} \; ; \quad 
V_r(\xi):= \frac{\xi\,(\sh\xi+\sin\xi) }{ 2\,(\ch\xi - \cos\xi)} \; .
\end{equation*}
Using these auxiliary functions, (\ref{f.exprq'}) reads\,: 
\begin{equation*}  
q_s(w,\beta , x , \zeta , z) 
\end{equation*}
\begin{equation*}  
=\, \frac{1}{8\pi^3 s^6} \int_{0}^\infty \Re\!\left[ F(y)\, e^{\rt1\! \left[\frac{B^2_s}{s^3}\,U_i(\xi) + \frac{B'_s-x}{s^2} + \frac{\ss w^2+\beta^2}{2s}\, V_i(\xi)\!\right] y}\right] e^{- \frac{B^2_s}{s^3}\,U_r(\xi) -\frac{\ss w^2+\beta^2}{2s}\, V_r(\xi) } \, y^2\, dy      
\end{equation*}
\begin{equation*}  
= {\ts\frac{\tilde q_s(w,\beta , \zeta , z)}{(3 \pi/2)\, s^2}} \int_{0}^\infty \Re\!\left[ F(y)\, e^{\rt1\!\! \left[\frac{B^2_s}{s^3}\,U_i(\xi) + \frac{B'_s-x}{s^2} + \frac{\ss w^2+\beta^2}{2s}\, V_i(\xi)\!\right] 2 \xi^2}\right]\! e^{- \frac{B^2_s}{s^3}[U_r(\xi)-12] -\frac{\ss w^2+\beta^2}{2s}[V_r(\xi)-1] } \, \xi^5 d\xi    
\end{equation*}
\begin{equation}  \label{f.exprqs} 
=\, \frac{\tilde q_s(w,\beta , \zeta , z)}{3\pi\, s^2} \int_{0}^\infty G_s(\xi)\,  e^{- \frac{B^2_s}{s^3}[U_r(\xi)-12] -\frac{\ss w^2+\beta^2}{2s}[V_r(\xi)-1]} \, \xi^5\, d\xi\,  , 
\end{equation}
with 
\begin{equation*}  
G_s(\xi) := \, \big(F_r(\xi) + \rt1 F_i(\xi)\big)\, e^{\rt1 \Lambda_s(\xi)} +  \big(F_r(\xi) - \rt1 F_i(\xi)\big)\, e^{-\rt1 \Lambda_s(\xi)}
\end{equation*}
\begin{equation*}  
= \, 2 F_r(\xi) \cos\!\big(\Lambda_s(\xi)\big) - 2 F_i(\xi) \sin\!\big(\Lambda_s(\xi)\big) \, , 
\end{equation*}
where we set \quad $\Lambda_s(\xi)\, :=\, 2\,\xi^2\big({\ts \frac{B^2_s}{s^3}\,U_i(\xi) + \frac{B'_s-x}{s^2} + \frac{\ss w^2+\beta^2}{2s}\, V_i(\xi)}\big)$. \parn 
Note that the functions $\,F_r,F_i,U_r,U_i,V_r,V_i, G_s\,$ are all even.   

\section{Small time asymptotics for the density $\,q_s(w,\beta , x , \zeta , z)$} \label{s.asbdqs}  \indf 
   We shall here use the expression (\ref{f.exprqs}) for $\,q_s(w,\beta , x , \zeta , z)$, to derive its asymptotics as $\,s\sea 0$, proceeding by adapting  the saddle-point method (see [Co] for example). This will require the following asymptotics for the auxiliary function entering that expression. \par 

\subsection{Auxiliary asymptotics} \label{s.auxasy}

As $\,\xi\to \infty\,$ we have \quad  $ U_r(\xi) = 2\xi + 4 + \frac{4}{\xi} + \O(\frac{1}{\xi^3})\, ; \; U_i(\xi) = 1 +\O(\frac{1}{\xi})\, ; \; V_r(\xi) = \frac{\xi}{2}[{\ss1+\O(e^{-\xi})}] \, ; \;$  $V_i(\xi) = \frac{1+\O(e^{-\xi})}{4\xi}$ ; \  $F(2\xi^2) = \O\big(\frac{e^{-\xi}}{\xi}\big)$\raise0.7pt\hbox{.} \parsn
Then near 0 we successively have\,: 
\begin{equation*}
U_r(\xi) = \frac{2\,\xi^3 \times \frac{\xi^3}{6}\big(1-\frac{17}{420}\xi^4+\O(\xi^8)\big) }{ \big(\frac{\xi^5}{30}\big)^2 + \big(\frac{\xi^3}{6}\big(1-\frac{17}{420}\xi^4\big)\big)^2 +\O(\xi^{10})} = \frac{\frac{\xi^6}{3}\big(1-\frac{17}{420}\xi^4+\O(\xi^8)\big) }{ \frac{\xi^6}{36}\big(1-\frac{17}{210}\xi^4\big) +\O(\xi^{10})} = 12 + {\ts \frac{17}{35}}\,\xi^4+\O(\xi^8) \, ;
\end{equation*}
\begin{equation*}
U_i(\xi) = \frac{\frac{\xi^6}{30} \big(1-\frac{113}{648}\xi^4+\O(\xi^8)\big)}{ \frac{\xi^6}{36}\big(1-\frac{17}{210}\xi^4+ \O(\xi^8)\big)} = {\ts \frac{6}{5}} - {\ts \frac{2119}{18900}}\,\xi^4+\O(\xi^8) \,  ; 
\end{equation*}
\begin{equation*}
V_i(\xi) =  {\ts \frac{1}{12}} - {\ts \frac{\xi^4}{756}} +\O(\xi^8) \, ; \  V_r(\xi) =  1 + {\ts \frac{\xi^4}{180}} +\O(\xi^8) \, ; 
\end{equation*}
\begin{equation*}  
F_r(\xi) = { \frac{6}{\xi^4}} - {\ts\frac{620659}{135600}} + \O(\xi^4) \,  ;  \   F_i(\xi) = \frac{1}{5} + \O(\xi^4)\, ;   
\end{equation*}
\begin{equation*}  
\Lambda_s(\xi)\, =\, 2  \left[\frac{B^2_s}{s^3}\big[{\ts \frac{6}{5}} - {\ts \frac{2119}{18900}}\,\xi^4+\O(\xi^8)\big] + \frac{B'_s-x}{s^2} + \frac{\ss w^2+\beta^2}{2s}\big[{\ts \frac{1}{12}} - {\ts \frac{1}{756}}\,\xi^4+\O(\xi^8) \big]\!\right] \xi^2 
\end{equation*}
\begin{equation*}  
= \, 2 \mu_s\, \xi^2 - \left[{\ts \frac{2119}{9450}\frac{B^2_s}{s^3} + \frac{1}{63}\frac{w^2+\beta^2}{12\,s}\,  }\right]\! \big[\xi^6 +\O(\xi^{10})\big] 
= \, 2 \mu_s\, \xi^2 \big[1 +\O(\xi^{4})\big],
\end{equation*}
where we have set  
\begin{equation}  \label{f.exprmus} 
\mu_s\, :=\, \frac{6B^2_s}{5\,s^3} + \frac{B'_s-x}{s^2} + \frac{ w^2+\beta^2}{24\,s}\, \raise1.5pt\hbox{.} 
\quad \hbox{ \  Let also }\quad   \nu_s := \frac{17\,B^2_s}{35\, s^3} + \frac{ w^2+\beta^2}{360\,s}\,\raise1.5pt\hbox{,} 
\end{equation}
\  so that \parn 
\begin{equation*}
e^{- \frac{B^2_s}{s^3}[U_r(\xi)-12] -\frac{\ss w^2+\beta^2}{2s}[V_r(\xi)-1]} = e^{
-\, \frac{B^2_s}{s^3}\big[{\ts \frac{17}{35}}\,\xi^4+\O(\xi^8)\big] -\frac{\ss w^2+\beta^2}{2s}\big[ {\ts \frac{1}{180}}\,\xi^4+\O(\xi^8)\big] } = e^{-\nu_s\,[\xi^4 +\O(\xi^8)]} . 
\end{equation*}
It is not difficult to see that actually \   $V_r(\xi) \ge 1\,$ and $\,U_r(\xi) \ge 12\,$,  for any real $\,\xi\,$. \parsn 
\if {\small  
Moreover, the signs of $\,[V_r(\xi) - 1]\,$ and $\,[U_r(\xi) - 12]\,$ are given respectively by\,: 
\begin{equation*}
2\,(\ch\xi - \cos\xi)\, [V_r(\xi) -1] =  \xi\,(\sh\xi+\sin\xi) - 2\,(\ch\xi - \cos\xi) \ge 0\, 
\end{equation*}
\big(since the  first derivatives are $\,\xi (\ch\xi + \cos\xi) - (\sh\xi+\sin\xi) $ and $\, \xi\,(\sh\xi-\sin\xi)\ge 0\,$\big), \parn
and  
\begin{equation*}
\5 (\ch\xi + \cos\xi){\ts\left[\left(\xi - \frac{\sh\xi+\sin\xi}{\ch\xi + \cos\xi}\right)^2 + \left(\frac{\sh\xi-\sin\xi}{\ch\xi + \cos\xi}\right)^2\right]} [U_r(\xi) - 12] 
\end{equation*}
\begin{equation*}  
=\, {\ts \xi^3 \,(\sh\xi-\sin\xi) - 6\, (\ch\xi + \cos\xi)\left[\left(\xi - \frac{\sh\xi+\sin\xi}{\ch\xi + \cos\xi}\right)^2 + \left(\frac{\sh\xi-\sin\xi}{\ch\xi + \cos\xi}\right)^2\right]}
\end{equation*}
\begin{equation*}  
=\, (\sh\xi-\sin\xi)\, \xi^3  - 6\,(\ch\xi + \cos\xi)\, \xi^2 + 12\, (\sh\xi+\sin\xi)\, \xi  - 12\, (\ch\xi-\cos\xi)\, \ge 0\,  , 
\end{equation*}
since its derivative is\,:
\begin{equation*}  
\big[(\ch\xi-\cos\xi)\, \xi - 3\,(\sh\xi-\sin\xi)\big]\xi^2 \ge 0\, , 
\end{equation*}
because (by the preceding)   
$$\frac{d}{d\xi}\big[(\ch\xi-\cos\xi)\, \xi - 3\,(\sh\xi-\sin\xi)\big] = \xi\,(\sh\xi+\sin\xi) - 2\,(\ch\xi - \cos\xi) \ge 0\, .  $$ 
}\fi    
Furthermore, 
\begin{equation*}  
U_r(\xi) -2\xi -4 = \frac{4\xi \big(1- \frac{2(e^{-\xi}-\sin\xi)}{\ch\xi + \cos\xi}\big) - 8 \big(\frac{\ch\xi - \cos\xi}{\ch\xi + \cos\xi}\big) -2\,\xi^3 \big(\frac{e^{-\xi}+ \cos\xi +\sin\xi}{\ch\xi + \cos\xi}\big) -4\,\xi^2 \big(\frac{e^{-\xi}+ \cos\xi -\sin\xi}{\ch\xi + \cos\xi}\big) }{ (\xi - 1)^2 +1 + 2\xi \big(\frac{e^{-\xi}+ \cos\xi -\sin\xi}{\ch\xi + \cos\xi}\big) - \frac{4\, \cos\xi }{\ch\xi + \cos\xi}} 
\end{equation*}
\begin{equation*}  \small \ts 
\left( =\, \frac{4\xi - 8 +\O(\xi^3e^{-\xi})  }{ \xi^2 -2\xi + 2 + \O(\xi\,e^{-\xi}) } 
\right)
\end{equation*}
\begin{equation*}
> \frac{4\xi - 8  - \frac{ 4 (\xi+1)^3}{\ch\xi + \cos\xi} }{ (\xi - 1)^2 +1 + \frac{ 4 (\xi+1)}{\ch\xi + \cos\xi} } > \frac{4\xi - 8  - 9 (\xi+1)^3 e^{-\xi} }{ (\xi - 1)^2 +1 +  9 (\xi+1) e^{-\xi} } > \frac{4\xi - 9 }{  \xi^2 -2\xi + 2.2 } \ge \frac{2}{  \xi} \times \frac{2\pi - 9 }{\pi -1 +  \frac{11}{20\pi} } 
\end{equation*}
for $\,\xi\ge 2\pi$, \  and similarly 
\begin{equation*}
V_r(\xi) - \frac{\xi}{ 2} = \frac{\xi\,(\cos\xi+\sin\xi - e^{-\xi}) }{ 2\,(\ch\xi - \cos\xi)} >  \frac{-\xi}{ \ch\xi - \cos\xi} > -3\xi\, e^{-\xi} > \frac{-1}{20}\,\raise1.7pt\hbox{,}
\end{equation*}
so that \  $V_r(\xi) > \pi - \frac{1}{20} >3\, $ for $\,\xi\ge 2\pi$. \parn
Therefore, in a small neighbourhood $\,[2\pi-\e ,\infty[ + \rt1 [-\e ,\e ]\,$ of $\,[2\pi,\infty[ \,$ we shall have 
\begin{equation}  \label{f.estUVr}
\Re\big[U_r(\xi) -2\xi \big] \ge 4 \quad \hbox{ and } \quad \Re\big[V_r(\xi) \big] \ge 3\, . 
\end{equation}
\vspace{-6mm} 

\subsection{Changes of contour and saddle-point method} \label{s.auxasy} \indf 
   Note that by Lemmas \ref{lem.zthz}, \ref{lem.explPsi}, (\ref{f.Psi-lambd}) and the changes of variable\,: $\,\lambda = \rt1 y = \pm 2\rt1 \xi^2$, the poles of the integrand in (\ref{f.exprqs}) are located at $\,e^{\rt1\! k\pi/2} (1+\rt1\!)\, y_n\,$ and $\,e^{\rt1\! k\pi/2} (1+\rt1\!)(n+1)\,\pi $, with $\,n, k\in\N$ and $\,\frac{17\pi}{12} < y_0 < y_1 < \ldots\,$, so that this integrand  is analytic at 0 with convergence radius 
$\sqrt{2}\,\pi$\raise0pt\hbox{.} \parn
In particular, we may perform the following changes of contour in (\ref{f.exprqs})\,: 
\begin{equation*}
3\pi\, s^2\, \frac{q_s(w,\beta , x , \zeta , z)}{\tilde q_s(w,\beta , \zeta , z)}\, = \int_{0}^\infty G_s(\xi)\,  e^{- \frac{B^2_s}{s^3}[U_r(\xi)-12] -\frac{\ss w^2+\beta^2}{2s}[V_r(\xi)-1]} \, \xi^5\, d\xi 
\end{equation*}
\begin{equation*}  
= \int_{0}^{(1+\rt1\!)\mu_s^{-\eta}}\! \big(F_r(\xi) + \rt1 F_i(\xi)\big)\, e^{\rt1 \Lambda_s(\xi)- \frac{B^2_s}{s^3}[U_r(\xi)-12] -\frac{\ss w^2+\beta^2}{2s}[V_r(\xi)-1]} \, \xi^5\, d\xi 
\end{equation*}
\begin{equation*}  
+ \int_{0}^{(1-\rt1\!)\mu_s^{-\eta}}\! \big(F_r(\xi) - \rt1 F_i(\xi)\big)\, e^{-\rt1 \Lambda_s(\xi)- \frac{B^2_s}{s^3}[U_r(\xi)-12] -\frac{\ss w^2+\beta^2}{2s}[V_r(\xi)-1]} \, \xi^5\, d\xi 
\end{equation*}
\begin{equation*}  
+ \int_{(1+\rt1\!)\mu_s^{-\eta}}^{(1+\rt1\!)\mu_s^{-\eta}+2\pi} \! \big(F_r(\xi) + \rt1 F_i(\xi)\big)\, e^{\rt1 \Lambda_s(\xi)- \frac{B^2_s}{s^3}[U_r(\xi)-12] -\frac{\ss w^2+\beta^2}{2s}[V_r(\xi)-1]} \, \xi^5\, d\xi 
\end{equation*}
\begin{equation*}  
+ \int_{(1-\rt1\!)\mu_s^{-\eta}}^{(1-\rt1\!)\mu_s^{-\eta}+2\pi}\! \big(F_r(\xi) - \rt1 F_i(\xi)\big)\, e^{-\rt1 \Lambda_s(\xi)- \frac{B^2_s}{s^3}[U_r(\xi)-12] -\frac{\ss w^2+\beta^2}{2s}[V_r(\xi)-1]} \, \xi^5\, d\xi 
\end{equation*}
\begin{equation*}  
+ \int_{(1+\rt1\!)\mu_s^{-\eta}+2\pi}^{(1+\rt1\!)\mu_s^{-\eta}+\infty} \! \big(F_r(\xi) + \rt1 F_i(\xi)\big)\, e^{\rt1 \Lambda_s(\xi)- \frac{B^2_s}{s^3}[U_r(\xi)-12] -\frac{\ss w^2+\beta^2}{2s}[V_r(\xi)-1]} \, \xi^5\, d\xi 
\end{equation*}
\begin{equation*}  
+ \int_{(1-\rt1\!)\mu_s^{-\eta}+2\pi}^{(1-\rt1\!)\mu_s^{-\eta}+\infty}\! \big(F_r(\xi) - \rt1 F_i(\xi)\big)\, e^{-\rt1 \Lambda_s(\xi)- \frac{B^2_s}{s^3}[U_r(\xi)-12] -\frac{\ss w^2+\beta^2}{2s}[V_r(\xi)-1]} \, \xi^5\, d\xi 
\end{equation*}
\begin{equation*}
=:\, J^0_s + \bar J^0_s +  J^\pi_s + \bar J^\pi_s + J^\infty_s + \bar J^\infty_s\, , 
\end{equation*}
where $\,\eta >\frac{1}{4}$ will be specified further. \parn
Note that the estimate (\ref{f.estUVr}) and the control $F(2\xi^2) = \O\big(\frac{e^{-\xi}}{\xi}\big)$   ensure the vanishing of the unmentioned limiting contribution (for any large enough fixed $s$, provided $\,\lim_{s\to 0}\limits \mu_s = \infty$) in the above changes of contour\,: 
\begin{equation*}  
\lim_{R\to\infty} \int_{R}^{R\pm (1+\rt1\!)\mu_s^{-\eta}}\! \big(F_r(\xi) + \rt1 F_i(\xi)\big)\, e^{\rt1 \Lambda_s(\xi)- \frac{B^2_s}{s^3}[U_r(\xi)-12] -\frac{\ss w^2+\beta^2}{2s}[V_r(\xi)-1]} \, \xi^5\, d\xi \, = 0\, .
\end{equation*}
   
   Now on the one hand, setting $\,\xi = (1\!+\!\rt1\!)\mu_s^{-\eta}\, t\,$ by the above we have\,:\parsn 
\begin{equation*}
J^0_s = {\frac{8}{\mu_s^{6\eta}}}\!\int_{0}^{1} \big[F_i(\xi) - \rt1\! F_r(\xi)\big]\, e^{\big[\rt1\! \Lambda_s(\xi) - \frac{B^2_s}{s^3}[U_r(\xi)-12] -\frac{\ss w^2+\beta^2}{2s}[V_r(\xi)-1]\big] } \,t^5 dt
\end{equation*}
\begin{equation*}
= {\frac{8}{5\mu_s^{6\eta}}}\!\int_{0}^{1} \Big[1+ {\ts\frac{15\rt1}{2\mu_s^{-4\eta}t^4}} + {\ts\frac{620659\rt1}{27120}} + \O(\mu_s^{-4\eta})\Big]\, 
e^{-4\big(\mu_s^{1-2\eta}\, t^2 - \nu_s \mu_s^{-4\eta}\big)\big(1+\O(\mu_s^{-4\eta})\big)} \,t^5\, dt
\end{equation*}
\begin{equation*}
= {\frac{8}{5\mu_s^{6\eta}}} \big(1+(\mu_s^{1-2\eta}+\nu_s)\O(\mu_s^{-4\eta})\big)\! \int_{0}^{1} \Big[1+ {\ts\frac{15\rt1}{2\mu_s^{-4\eta}t^4}} + {\ts\frac{620659\rt1}{27120}} + \O(\mu_s^{-4\eta})\Big] 
e^{-4\mu_s^{1-2\eta} t^2} \,t^5\, dt\, , 
\end{equation*}
so that for \  $\frac{1}{4} < \eta <\frac{1}{2}\,$: 
\begin{equation*}
J^0_s + \bar J^0_s\, =\, {\frac{16}{5\mu_s^{6\eta}}} \big(1+\O(\nu_s\mu_s^{-4\eta})\big)\int_{0}^{1}  
e^{-4\mu_s^{1-2\eta}\, t^2} \,t^5\, dt \, 
\end{equation*}
\begin{equation*}
= \, {\frac{1+\O(\nu_s\,\mu_s^{-4\eta})}{40\,\mu_s^{3}}} \int_{0}^{4\mu_s^{1-2\eta}}\! e^{- u} \,u^2 du \, 
= \, {\frac{1+\O(\mu_s^{1-4\eta})}{20\,\mu_s^{3}}}\, \raise1.9pt\hbox{,} 
\end{equation*}
provided  \  $\,\lim_{s\to 0}\limits\, \mu_s = \infty\,$ and $\,\nu_s = \O(\mu_s)$. \pars 

   Then we use the saddle-point method (as described in [Co]) to deal with the intermediate part $(J^\pi_s + \bar J^\pi_s)$.  \  By the above we know that we have 
\begin{equation*}
U_i(\xi) = {\ts \frac{6}{5}} + \xi^4\, \tilde U_i(\xi)\, , \; V_i(\xi) = {\ts \frac{1}{12}} + \xi^4\, \tilde V_i(\xi)\, , \;  U_r(\xi)-12 = \xi^4\, \tilde U_r(\xi)\, , \; V_r(\xi)-1 = \xi^4\, \tilde V_r(\xi)\, , \; 
\end{equation*}
with even functions $\,\tilde U_i\, , \tilde V_i\, , \tilde U_r\, , \tilde V_r\,$ that are analytic outside the sequence \parn
$\big\{e^{\rt1\, k\pi/2} (1+\rt1\!)\, y_n\,\big|\,n, k\in\N, \, \frac{17\pi}{12} < y_0 < y_1 <\ldots \big\}$ and then converge in the compact disc (centred at 0) of radius $\frac{17\pi}{6\sqrt{2}}> 2\pi $\raise0pt\hbox{.} \big(Note that the proof of Lemma \ref{lem.explPsi} shows that the points $\,e^{\rt1\! k\pi/2} (1+\rt1\!)(n+1)\,\pi \,$ are poles only for the functions $F_r,F_i\,$.\big) Therefore we can write the phase as follows\,: 
\begin{equation*}  
\f(\xi) := \, {\rt1 \Lambda_s(\xi) - \frac{\ss B^2_s}{s^3}[U_r(\xi)-12] -\frac{\ss w^2+\beta^2}{2s}[V_r(\xi)-1]} 
\end{equation*}
\begin{equation*}
=\, 2\rt1\mu_s\, \xi^2  - \big( {\ts\frac{B^2_s}{s^3}\, [\tilde U_r(\xi)-2\rt1 \tilde U_i(\xi)] + \frac{\ss w^2+\beta^2}{2s}\, [\tilde V_r(\xi)}-2\rt1 \tilde V_i(\xi)]\big)\, \xi^4\, . 
\end{equation*}
Setting \  $M:= \max\big\{ |w(\xi)|\,\big|\,\xi\in\C,  |\xi| = \frac{17\pi}{6\sqrt{2}}\big\}$, by Cauchy's inequality, for $\,|\xi|\le \pi+\frac{17\pi}{12\sqrt{2}}\,$ we have 
\begin{equation*}
\big| \f(\xi) - 2\rt1\mu_s\, \xi^2\big|\,  \le \, \frac{M\, |\xi|^4}{\frac{289\pi^2}{72}(\frac{289\pi^2}{72}-|\xi|^2)} \,  \le \, \frac{4\, (72)^2\,M\, |\xi|^4}{289(579-408\sqrt{2})\,\pi^4} \,<\,1156\,M . 
\end{equation*}
As before near 0 we have \  $\big(F_r(\xi) + \rt1 F_i(\xi)\big) \xi^4  = 6+ \O(\xi^4)$. \quad  Hence for small $\,s\,$ 
\begin{equation*}  
J^\pi_s\, = \int_{(1+\rt1\!)\mu_s^{-\eta}}^{(1+\rt1\!)\mu_s^{-\eta}+2\pi} \! \big(F_r(\xi) + \rt1 F_i(\xi)\big)\, e^{\rt1 \Lambda_s(\xi)- \frac{B^2_s}{s^3}[U_r(\xi)-12] -\frac{\ss w^2+\beta^2}{2s}[V_r(\xi)-1]} \, \xi^5\, d\xi  
\end{equation*}
\begin{equation*}
= \int_{0}^{\pi} \Big[6 + \O\big(|(1\!+\!\rt1\!)\mu_s^{-\eta} + t|^4\big)\Big]\, e^{2\rt1\mu_s \big((1\!+\!\rt1\!)\mu_s^{-\eta} +t\big)^2 + \O(1)} \big((1\!+\!\rt1\!)\mu_s^{-\eta} + t\big) dt
\end{equation*}
\begin{equation*}
= \,\O(1) \int_{0}^{\pi} e^{-4\mu_s^{1-2\eta} -4\mu_s^{1-\eta}t} \, dt\, = \,\O\big(e^{-4\mu_s^{1-2\eta}}\big/\mu_s^{1-\eta}\big)\, . 
\end{equation*}
The same of course holds for  $\,\bar J^\pi_s$.   \par  
To deal with the remaining contribution $(J^\infty_s + \bar J^\infty_s)$, we use the lower estimates near infinity\,: (\ref{f.estUVr}) computed above, as follows\,:
\begin{equation*}  
J^\infty_s\, = \int_{(1+\rt1\!)\mu_s^{-\eta}+2\pi}^{(1+\rt1\!)\mu_s^{-\eta}+\infty} \! \big(F_r(\xi) + \rt1 F_i(\xi)\big)\, e^{\rt1 \Lambda_s(\xi)- \frac{B^2_s}{s^3}[U_r(\xi)-12] -\frac{\ss w^2+\beta^2}{2s}[V_r(\xi)-1]} \, \xi^5\, d\xi  
\end{equation*}
\begin{equation*}
= \int_{2\pi}^{\infty}  \O\bigg[{\frac{e^{- t}}{ t}}\bigg] \, e^{ - \frac{B^2_s}{s^3}[2(\mu_s^{-\eta} +t) - 8] -\frac{\ss w^2+\beta^2}{s}} (1 + t)^5 dt\,   =\, \O(1)\, e^{\frac{8 B^2_s}{s^3}-\frac{\ss w^2+\beta^2}{s}} \int_{2\pi}^{\infty}   e^{- \big[\frac{2 B^2_s}{s^3} + 1\big] t}\, t^4\, dt
\end{equation*}
\begin{equation*}
=\, \O(1)\, e^{\frac{8 B^2_s}{s^3}-\frac{\ss w^2+\beta^2}{s}} \big({\ts\frac{2 B^2_s}{s^3} + 1}\big)^{-5} \int_{2\pi \big[\frac{2 B^2_s}{s^3} + 1\big]}^{\infty}   e^{- t}\, t^4\, dt\,  =\,  \O\Big[\big({\ts\frac{2 B^2_s}{s^3} + 1}\big)^{-1}\Big]\, e^{-4(\pi-2)\frac{B^2_s}{s^3}-\frac{\ss w^2+\beta^2}{s}} , 
\end{equation*}
and the same of course holds for  $\,\bar J^\infty_s$. \parm 
   So far, we have obtained\,:   
\begin{equation*}
3\pi\, s^2\, \frac{q_s(w,\beta , x , \zeta , z)}{\tilde q_s(w,\beta , \zeta , z)}\, =\, J^0_s + \bar J^0_s +  J^\pi_s + \bar J^\pi_s + J^\infty_s + \bar J^\infty_s
\end{equation*}
\begin{equation*}
= \, {\frac{1+\O(\nu_s\,\mu_s^{-4\eta})}{20\,\mu_s^{3}}} + \O\big(e^{-4\mu_s^{1-2\eta}}\big/\mu_s^{1-\eta}\big) + \O\Big[\big({\ts\frac{2 B^2_s}{s^3} + 1}\big)^{-1}\Big]\, e^{-4(\pi-2)\frac{B^2_s}{s^3}-\frac{\ss w^2+\beta^2}{s}}
\end{equation*}
\begin{equation*}
= \, {\frac{1+\O(\mu_s^{1-4\eta})}{20\,\mu_s^{3}}} \quad \hbox{ for }\; {\ts\frac{1}{4} < \eta <\frac{1}{2}}\, \raise1.3pt\hbox{,} 
\end{equation*}
provided both \  $\,\lim_{s\to 0}\limits\, \mu_s = \infty\,$ and $\,\nu_s = \O(\mu_s)$. \parn
By (\ref{f.defB2B'}) and (\ref{f.exprmus}), this condition  holds as soon as  
both \  $\lim_{s\to 0}\limits \mu_s = \infty\,$ and  $\,\frac{x}{s^2} \le \frac{z^2+\zeta^2}{s^3} + \frac{\mu_s}{\e}$ (for some $\e>0$), and then also as soon as both \  
$\,\lim_{s\to 0}\limits \nu_s = \infty\,$ and $\,2s\, x  \le  (z^2+\zeta^2)+\e s^2(w^2+\beta^2) $. \parsn
Finally, under this condition, as $\,s\sea 0\,$  and for any positive $\e\,$, we have obtained\,:
\begin{equation*}
q_s(w,\beta , x , \zeta , z)\, =\, \frac{1+\O\big(\mu_s^{\e-1}\big)}{60\,\pi\,s^2\,\mu_s^{3}}\times \tilde q_s(w,\beta , \zeta , z)\, . 
\end{equation*}
The result of this section (and main result) is thus the following off-diagonal equivalent. 
\bThe  \label{pro.asymptqs} \  As  $\,s\sea 0$, for any positive $\e$, uniformly for $\,x\ge 0\,$ and $(w,\beta , \zeta , z)\in\R^4$ such that\,:
\begin{equation}  \label{f.condit} 
\mu_s \, \equiv\, 
\frac{3[( z-\frac{sw}{12})^2 + (\zeta-\frac{s\beta}{12})^2]}{5\, s^3} + \frac{ w^2+\beta^2 }{16\, s} - \frac{x}{s^2} \lra \infty\quad \hbox{and} \quad  \frac{x}{s^2} \le \frac{z^2+\zeta^2}{s^3} + \frac{\mu_s}{\e}\, \raise1.5pt\hbox{,}
\end{equation}
we have  
\begin{equation*}
q_s(w,\beta , x , \zeta , z)\, =\, \frac{1+\O\big(\mu_s^{\e-1}\big)}{20\,\pi^3\,s^6\,\mu_s^{3}}\times  \exp\left[- { \frac{3  \left[(sw-2z)^2+(s\beta-2\zeta)^2\right] }{2\,s^3}} - {\frac{w^2+\beta^2}{2\, s}} \right] .  
\end{equation*}
An alternative condition (to (\ref{f.condit}) above) guaranteeing this asymptotic equivalent is
\begin{equation}  \label{f.condit'} 
\frac{(sw-2z)^2+(s\beta-2\zeta)^2}{s^3} + \frac{w^2+\beta^2}{s}\, \lra \infty\quad \hbox{and} \quad  2s\, x  \le  (z^2+\zeta^2)+ \e s^2(w^2+\beta^2)\, . 
\end{equation}
\eThe
\bRem  \label{rem.finale} \  {\rm  
   Theorem \ref{pro.asymptqs} rather precisely yields the small time asymptotic behaviour of the heat kernel of the second chaos approximation $(Y_s)$ to the Dudley relativistic diffusion $(X_s)$. But this is not enough to derive any small time asymptotics for the density $\,p_s(\lambda, \mu, x , y , z)$ of the original process $\,X_s\equiv (\lambda_s\,,\, \mu_s\,,\, x_s\,,\, y_s\,,\, z_s)$, even for fixed $(\lambda, \mu, x , y , z)$. The reason is that the computations performed in Section \ref{sec.FLTransf} above cannot work beyond the second chaos, so that Section \ref{sec.sta} cannot yield a sufficiently precise control on the gap between both tangent processes $(X_s)$ and $(Y_s)$.  To be more specific, Section \ref{sec.sta} and Theorem \ref{pro.asymptqs} heuristically yield \,: 
\begin{equation*}
p_s(w,\beta , x , \zeta , z) \approx \P\big[ \lambda_s=w\,,\, \mu_s=\beta\,,\, x_s=x\,,\, y_s=\zeta\,,\, z_s=z\big]  \approx 
\end{equation*}
\begin{equation*}
\P\Big[ ( w_s\, , \, \beta_s\, , \,A_s\, , \,\zeta_s\, , \,\bar z_s) = (w,\beta , x , \zeta , z) + \O(R) (s^{3/2} , s^{3/2} , s^{5/2} , s^{5/2} , s^{5/2})\Big]  + \O\big(e^{-R^\kappa/c}\big)
\end{equation*}
\begin{equation*}
\approx \frac{\exp\!\left[ - { \frac{3  \left[\left(sw-2z+ \O(s^{5/2})\right)^2+\left(s\beta-2\zeta+ \O(s^{5/2})\right)^2\right] }{2\,s^3}} - {\frac{w^2+\beta^2 + \O(s^{3/2})}{2\, s}}\right]  }{ 20\,\pi^3\,s^6\,\mu_s\big(w+ \!\O(s^{3/2}) ,\beta+\! \O(s^{3/2}) , x+\! \O(s^{5/2}) , \zeta+\! \O(s^{5/2}) , z+\! \O(s^{5/2})\big)^{3}} + \O\big(e^{-R^\kappa/c}\big)
\end{equation*}
\begin{equation*}
=\, \exp\!\left[ - { \frac{3  \left[\left(sw-2z\right)^2+\left(s\beta-2\zeta\right)^2\right] }{2\,s^3}} - {\frac{w^2+\beta^2}{2\, s}} + \O(s^{-1/2})\right] + \O\big(e^{-R^\kappa/c}\big)\, , 
\end{equation*}
which were not too bad only if the additive term  $\O\big(e^{-R^\kappa/c}\big)$ were not there to ruin such estimate. Indeed, even taking $R\asymp s^{-\gamma}$ would only control this correction term by (at best) $e^{- s^{-2\gamma/3}}$, which would be significant only for $\gamma\ge 9/2$, so that the remaining information would then reduce to nothing. 
}\eRem    
\vspace{-5mm} 

\if{  
\section{Small time asymptotic density of the process $\,X_s$} \label{sec.AsympDensXs} \indf 
   We now use Section \ref{sec.sta} and Theorem \ref{pro.asymptqs}, to derive small time asymptotics relating to the density $\,p_s\,$ of the original process $\,X_s\equiv (\lambda_s\,,\, \mu_s\,,\, x_s\,,\, y_s\,,\, z_s)$. \parn
\bThe  \label{pro.asympX} \  As $\,s\sea 0\,$, for any $(\lambda,\mu , x , y , z)\in (\R^2)^*\times\R_+^*\times\R^2\,$ we have\,:  
\begin{equation*}    ??? 
\end{equation*}
uniformly on the complement of any neighbourhood of $\{x=0\}\cup\{\lambda=\mu=0\}$. 
\eThe  
\ub{Proof} \quad     ???

Let $F$ be a $C^1$ non-negative non-zero bounded function on $\R^2\times\R_+\times\R^2$ having bounded gradient, which vanishes in a neighbourhood of $\{x=0\}\cup\{w=\beta=0\}$.  \   
By Lemma \ref{lem.approxXY}   
we have\,: 
\begin{equation*}
\int_{\sR^5} F(\lambda,\mu,x,y,z)\, p_s(\lambda,\mu,s+x,y,z)\, \ch\lambda\, d\lambda\, d\mu\, dx\, dy\, dz \, = \, \E_0\Big[F(\lambda_s,\mu_s,x_s-s,y_s,z_s)\Big] 
\end{equation*}
\begin{equation*}
= \, \E_0\Big[F(w_s,\beta_s,A_s , \zeta_s, \bar z_s) + \min\big\{2\|F\| , \|\nabla F\| (R_s+R'_s)\big\} \Big] 
\end{equation*}
\begin{equation*}
= \int F(w,\beta,x, \zeta , z)\, q_s(w,\beta , x , \zeta , z)\, dw\,d\beta\, dx\, d\zeta\, dz\, + 5\,\|F\|\, e^{-R^\kappa /c} + 2\, \|\nabla F\|\, R\, s^{3/2} \, , 
\end{equation*}
for any large $R$ and small $\,s\,$.  Now taking \ $R= \Big(c\log\big(\frac{a\kappa}{cS}\big)\Big)^{1/\kappa} \Big(1- (1-\kappa) \frac{\log [c\log (\frac{a\kappa}{cS} ) ]}{\kappa^2\log (\frac{a\kappa}{cS} )} \Big)$, we find that \  
$\min_{R>0}\limits \{a\, e^{-R^\kappa/c} + RS\} \le \, S \Big(c\log\big(\frac{a\kappa}{cS}\big)\Big)^{1/\kappa}$  for any small $\,S$, so that setting $\,c':= c\,\log\!\big(\frac{5\kappa\,\|F\|}{2c\, \|\nabla F\|}\big)$, \   we can replace the above quantity 
$\,5\,\|F\|\, e^{-R^\kappa /c} + 2\, \|\nabla F\|\, R\, s^{3/2}\,$ by $\,C_s(F)$, which will below denote any quantity satisfying $\, 0\le C_s(F)\le\,  2\, \|\nabla F\| \big({\ts\frac{3c}{2}}\log({\ts\frac{1}{s}}) + c'\big)^{1/\kappa} s^{3/2}$. \parn 
Thus for small $s\,$ and any positive $\,u = u_s(F)\,$ we have\,:  
\begin{equation*}
\int_{\sR^5} F(\lambda,\mu,x,y,z)\, p_s(\lambda,\mu,s+x,y,z)\, \ch\lambda\, d\lambda\, d\mu\, dx\, dy\, dz  
\end{equation*}
\begin{equation*}
= \int F(w,\beta,x, \zeta , z)\, q_s(w,\beta , x , \zeta , z)\, dw\,d\beta\, dx\, d\zeta\, dz\, + C_s(F)\, = \int \big(F+C_s(F)\big)\, q_s\, d\Lambda   
\end{equation*}
\begin{equation*}
= \int_{\{F\ge u\}} F\times \big(1+ {\ts\frac{C_s(F)}{u}}\big)\, q_s\, d\Lambda  + \int_{\{F< u\}} \big(F+C_s(F)\big)\, q_s\, d\Lambda
\end{equation*}
\begin{equation*}
= \int F\times \big(1+ {\ts\frac{C_s(F)}{u_s(F)}}\big)\, q_s\, d\Lambda\,  + C_s(F)\! \int \big(1-{\ts\frac{F}{u_s(F)}}\big)^+\, q_s\, d\Lambda\, . 
\end{equation*}
So far, it remains to find $\,u = u_s(F)$ such that \big(majorising $(1-\frac{F}{u})^+$ by $\,1_{\{F<u\}}$\big), eventually as $s\sea 0$\,: 
\begin{equation*}
 \|\nabla F\|\,  s^{3/2-\e} \int_{\{F< u\}} \, q_s\, d\Lambda\, \le  \int F \, q_s\, d\Lambda\, .
\end{equation*}
Alternatively, restricting to compactly supported (?) non negative $C^1$ functions $F$ such that $\|\nabla F\| =1$, we have to find $\,v = v_s(F) \equiv \frac{u_s(F)}{\|\nabla F\|}\,$ such that 
\begin{equation*}
\limsup_{s\sea 0} \bigg[  s^{3/2-\e} \int_{\{F< v\}} \, q_s\, d\Lambda\,\bigg/ \int F \, q_s\, d\Lambda\bigg] \le  1\, . 
\end{equation*}
According to Theorem \ref{pro.asymptqs}, this amounts to finding the values of $\,v = v_s(F)$ such that\,: 
\begin{equation*}
1\,\ge \, \limsup_{s\sea 0}\, \frac{s^{3/2-\e} \int_{\{F< v\}} \, \mu_s^{-3} \,  e^{- { \frac{3  \left[(sw-2z)^2+(s\beta-2\zeta)^2\right] }{2\,s^3}} - {\frac{w^2+\beta^2}{2\, s}} } dw\,d\beta\, dx\, d\zeta\, dz }{\int F \, \mu_s^{-3}\,  e^{- { \frac{3  \left[(sw-2z)^2+(s\beta-2\zeta)^2\right] }{2\,s^3}} - {\frac{w^2+\beta^2}{2\, s}} } dw\,d\beta\, dx\, d\zeta\, dz}
\end{equation*}
\begin{equation*}
= \, \limsup_{s\sea 0}\, \frac{s^{3/2-\e} \int_{\big\{F\big(w\sqrt{s}\,,\beta\sqrt{s}\,,xs^2, \zeta\sqrt{s^3}\, , z\sqrt{s^3}\,\big)< v\big\}} \, \mu^{-3} \,  e^{- \frac{3}{2} \left[(w-2z)^2+(\beta-2\zeta)^2\right] - {\frac{w^2+\beta^2}{2}} }\, d\Lambda }{\int F\big(w\sqrt{s}\,,\beta\sqrt{s}\,,xs^2, \zeta\sqrt{s^3}\, , z\sqrt{s^3}\,\big) \mu^{-3}\,  e^{- \frac{3}{2} \left[(w-2z)^2+(\beta-2\zeta)^2\right] - {\frac{w^2+\beta^2}{2}} }\, d\Lambda}
\end{equation*}

where \quad  $ \mu \equiv\mu(w,\beta,x, \zeta , z) := \,\frac{3 ( z^2 + \zeta^2)}{5} - \frac{wz + \beta\zeta + 10\,x}{10} + \frac{w^2+\beta^2}{15}\,\raise1.5pt\hbox{,}$ 

$$  $$ 

Rque\, :  \quad 
${\ds \mu_s \, =\, \frac{3[( z-\frac{sw}{12})^2 + (\zeta-\frac{s\beta}{12})^2]}{5\, s^3} + \frac{ w^2+\beta^2 }{16\, s} - \frac{x}{s^2} \, }$ 

$$ $$ 

Mais  \   ATTENTION  \`a la zone o\`u le Theorem \ref{pro.asymptqs} ne dit rien ... !! , id est hors de 
\begin{equation*}   
\mu_s \, \equiv\,\frac{3 ( z^2 + \zeta^2)}{5\, s^3} - \frac{wz + \beta\zeta + 10\,x}{10\, s^2} + \frac{w^2+\beta^2}{15\, s}\, \lra \infty\quad \hbox{and} \quad  \frac{x}{s^2} \le \frac{z^2+\zeta^2}{s^3} + \frac{\mu_s}{\e}\, 
\end{equation*}
ou bien hors de 
\begin{equation*}   
\frac{(sw-2z)^2+(s\beta-2\zeta)^2}{s^3} + \frac{w^2+\beta^2}{s}\, \lra \infty\quad \hbox{and} \quad  2s\, x  \le  (z^2+\zeta^2)+ \e s^2(w^2+\beta^2)\, . 
\end{equation*}

Cela s'arrange peut-\^etre en prenant $F$ \`a support compact... ou alors (mieux ?) 
$$(z^2+\zeta^2) \ge \e\,$$ 

\eject  
}\fi  


\section{Proofs of some technical results} \label{sec.Proofs} 
   We gather here the rather technical proofs of Lemmas \ref{lem.approxXY} and \ref{lem.LaplZ} and Proposition \ref{pro.abscCv}. 
\pars    \parmn    
\ub{\bf  Proof of Lemma \ref{lem.approxXY}}  \   $(i)$ \    Equation (\ref{f.sde}) entails that for small proper time $\,s\,$ we have$\,$: 
$$ \lambda_s = w_s + {\ts\frac{1}{2}}\! \int_0^s \th\!\big( w_\tau + o(\tau)\big) d\tau = w_s + {\ts\frac{1}{2}}\! \int_0^s w_\tau\,  d\tau + o(s^{2}) = \, w_s + o(s^{3/2-\e}) $$ 
$$= w_s + {\ts\frac{1}{2}}\! \int_0^s \th\!\big( w_\tau + o(\tau^{3/2-\e})\big) d\tau = w_s + {\ts\frac{1}{2}}\! \int_0^s\! \big( w_\tau + o(\tau^{3/2-\e})\big) d\tau =  w_s + {\ts\frac{1}{2}}\!\int_0^s\! w_\tau\,d\tau  + o(s^{5/2-\e}). $$ 
Then 
$$ \ch \lambda_s = \ch\!\big[ w_s + o(s^{3/2-\e})\big] = 1 + {\ts\frac{1}{2}}\,w_s^2 + o(s^{2-\e}) = 1 + o(s^{1-\e}) \, ,  $$ 
and by Equation (\ref{f.sde'}) we have$\,$: 
$$ \dot z_s = \sh \lambda_s = \sh\!\big[ w_s + o(s^{3/2-\e})\big] = w_s  + o(s^{3/2-\e}) \, .   $$ 
\begin{equation*}  \small \ts 
\left(\hbox{Also, } \quad \dot z_s = \int_0^s \sqrt{1+ \dot z_\tau ^2}\; dw_\tau +  \int_0^s\dot z_\tau\, d\tau\, . \right) 
\end{equation*}
Similarly,  
$$ \mu_s = \int_0^s \frac{d\beta_\tau}{\ch[w_\tau  + o(\tau^{3/2-\e})]} =  \int_0^s \big(1- o(\tau^{1-\e})\big) d\beta_\tau = \beta_s + o(s^{3/2-\e})\, , $$ 
$$ \ch \mu_s =  1 + o(s^{1-\e}) \, ,  \quad  \sh \mu_s =  \beta_s  + o(s^{3/2-\e}) \, .   $$ 
\begin{equation*}  \small \ts 
\left(\hbox{Also, } \quad \sh \mu_s = \int_0^s {\ts\sqrt{\frac{1+ \sh^2\mu_\tau}{1+ \dot z_\tau ^2}}}\; d\beta_\tau +  \int_0^s {\ts\frac{\sh\mu_\tau}{2(1+ \dot z_\tau^2)}}\, d\tau\, . \right) 
\end{equation*}
Finally, the result follows at once from : 
$$ x_s = \int_0^s \big[ 1 + {\ts\frac{1}{2}}\,w_\tau^2 + o(\tau^{2-\e})\big] \big[ 1 + {\ts\frac{1}{2}}\,\beta_\tau^2 + o(\tau^{2-\e})\big] d\tau = s + {\ts\frac{1}{2}} \int_0^s \big[\beta_\tau^2 + w_\tau^2\big] d\tau + o(s^{3-\e})\, ; $$ 
$$ y_s = \int_0^s \big[ 1 + o(\tau^{1-\e})\big] \big[ \beta_\tau + o(\tau^{3/2-\e}) \big] d\tau = \int_0^s\! \beta_\tau\,d\tau + o(s^{5/2-\e})\, ; \; z_s = \int_0^s\! w_\tau\,d\tau + o(s^{5/2-\e})\, . $$ 
\begin{equation*} \small \ts 
\left( \hbox{More precisely,}  \quad z_s = \int_0^s w_\tau\,d\tau + \int_0^s\!  \int_0^\tau\! w_u\, du \,d\tau +\5\int_0^s\! \int_0^\tau\! w_u^2\, dw_u \,d\tau + o(s^{7/2-\e})\, . \right) 
\end{equation*}

$(ii)$ \  For any $\,s\in [0,1]$ we \as have : \  
${\ds  |\lambda_s|\, \le\,  \sup_{[0,s]}\limits |w| + {\ts\frac{1}{2}}\! \int_0^s |\lambda_\tau|\, d\tau\,}$, whence by Gronwall's Lemma : ${\ds  \sup_{[0,s]}\limits |\lambda|\, \le\,  \sup_{[0,s]}\limits |w|\times e^{s/2}\,}$. Then for $\,0\le s\le 2\log 2\,$: \parn
\centerline{${\ds \Lambda_s :=  \sup_{[0,s]}\limits |\lambda- w|\, \le \5\int_0^s  \sup_{[0,\tau]}\limits |w|\times e^{\tau/2} d\tau\, \le \int_0^s  \sup_{[0,\tau]}\limits |w|\, d\tau \, . }$} \parn
Hence for $\,R\ge 1\,$ and $\,0\le s\le 2\log 2\,$: 
$$ \P\bigg[ \sup_{0\le t\le s}\! |\Lambda_t| \ge R\, s^{\frac{3}{2}}\bigg] \le \P\bigg[ \int_0^1\! \sup_{[0,\tau]}\limits |w|\, d\tau \ge R\bigg] \le \P\bigg[ \sup_{[0,1]}\limits |w| \ge R\bigg] \le 4\, \P[w_1\ge R] \le 2\, e^{-\frac{R^2}{2}} .   $$
Then \  ${\ds \beta_s - \mu_s = 2 \int_0^s \frac{\sh\!^2(\lambda_\tau/2)}{\ch \lambda_\tau}\, d\beta_\tau \,\equiv\, 2\, B\bigg[\int_0^s \frac{\sh\!^4(\lambda_\tau/2)}{\ch\!^2\lambda_\tau}\, d\tau\bigg]  }$, \   so that 
$$ \P\bigg[ \sup_{0\le t\le s}\! \big|\beta_t - \mu_t\big| \ge R\, s^{\frac{3}{2}}\bigg] \le 2\,\E\bigg[ \exp\bigg(-\frac{R^2 s^3}{8}\bigg/ \int_0^s \frac{\sh\!^4(\lambda_\tau/2)}{\ch\!^2\lambda_\tau}\, d\tau \bigg)  \bigg]  $$
$$ \le 2\,\E\bigg[ \exp\bigg(-2 R^2 s^3\bigg/ \int_0^s \lambda_\tau^4\, d\tau \bigg)  \bigg] 
 \le 2\,\E\bigg[ \exp\bigg(-\frac{R^2}{8\,\sup_{[0,1]}\limits |w|^4} \bigg)  \bigg]  $$ 
$$ = 2\int_0^1\P\bigg[  -\frac{R^2}{8\,\sup_{[0,1]}\limits |w|^4}  > \log y \bigg] dy 
= 2\int_0^1\P\bigg[ \sup_{[0,1]}\limits |w|^4  > \frac{-R^2}{8 \log y} \bigg] dy  $$ 
$$ \le  4\int_0^1e^{ -\5 \sqrt{\frac{-R^2}{8 \log y}}}\, dy  = 4\int_0^1e^{\frac{-R}{4  \sqrt{2 \log (1/y)}}}\, dy 
= 4\int_0^\infty e^{\frac{-R}{4  \sqrt{2 t}}-t }\, dt\,  \le\, 132 \, e^{-R^{2/3}/(32)^{1/3}} . $$ 
Then \qquad ${\ds \bigg| z_s - \int_0^s\! w_\tau\,d\tau\bigg| \, \le \int_0^s\Big| \sh\!\big(w_\tau + (\lambda_\tau -w_\tau)\big) - w_\tau\Big|\, d\tau }$  
$$ = \int_0^s\Big| \sh w_\tau - w_\tau + 2\,\sh\!\big[{\ts\frac{\lambda_\tau -w_\tau}{2}}\big]\,\ch\!\big[{\ts\frac{\lambda_\tau + w_\tau}{2}}\big]  \Big|\, d\tau\, \le  \int_0^s\Big( R^3 \tau^{3/2} + 2 \,\sh\!\big[{\ts\frac{R\, \tau^{3/2}}{2}}\big]\,\ch\!\big[2R\sqrt{\tau}\,\big]  \Big) d\tau\, $$ 
$$ \le  \int_0^s ( R^3 + 2R)\, \tau^{3/2} d\tau\, \le \, R^3 s^{5/2} \quad \hbox{ for }\, 0\le s\le s_R >0\, , \   \hbox{ with probability }\, 1-\O(e^{- {R^2/2}}) . $$ 
Similarly \qquad ${\ds \bigg| y_s - \int_0^s\! \beta_\tau\,d\tau\bigg| \, \le \int_0^s\Big| \ch\!\big(w_\tau + (\lambda_\tau -w_\tau)\big)\,\sh\!\big(\beta_\tau + (\mu_\tau-\beta_\tau) \big) - \beta_\tau\Big|\, d\tau }$  
$$ = \int_0^s\Big| \big(1+ 2\,\sh\!^2\big[{\ts\frac{w_\tau}{2}}\big] + 2\,\sh\!\big[{\ts\frac{\lambda_\tau -w_\tau}{2}}\big]\,\sh\!\big[{\ts\frac{\lambda_\tau + w_\tau}{2}}\big]\big) \big(\sh\beta_\tau + 2\,\sh\!\big[{\ts\frac{\mu_\tau - \beta_\tau}{2}}\big]\,\ch\!\big[{\ts\frac{\mu_\tau + \beta_\tau}{2}}\big]\big) - \beta_\tau \Big|\, d\tau $$ 
$$ \le  \int_0^s\Big( R^3 \tau^{3/2} + R\, \tau^{3/2} + \O\big(R\,\tau^{5/2}\big) \Big) d\tau\, \le \, R^3 s^{5/2} \quad \hbox{ for }\, 0\le s\le s_R' >0\, , $$ 
with probability $\, 1-\O(e^{- {R^{2/3}/4}})\,$. \quad 
Finally, in the same way we obtain\,:   
$$  \bigg| x_s - \Big( s+ {\ts\frac{1}{2}}\! \int_0^s\! \big[\beta_\tau^2 + w_\tau^2\big] d\tau\Big) \bigg| \le \int_0^s\Big| \ch\!\big(w_\tau + (\lambda_\tau -w_\tau)\big)\,\ch\!\big(\beta_\tau + (\mu_\tau-\beta_\tau) \big) -1 - \5  \big[\beta_\tau^2 + w_\tau^2\big] \Big|\, d\tau $$ 
$$ = \int_0^s\Big| \big(\ch w_\tau + 2\,\sh\!\big[{\ts\frac{\lambda_\tau -w_\tau}{2}}\big]\,\sh\!\big[{\ts\frac{\lambda_\tau + w_\tau}{2}}\big]\big) \big(\ch\beta_\tau + 2\,\sh\!\big[{\ts\frac{\mu_\tau - \beta_\tau}{2}}\big]\,\sh\!\big[{\ts\frac{\mu_\tau + \beta_\tau}{2}}\big]\big) -1 - \5\big[\beta_\tau^2 + w_\tau^2\big]\Big|\, d\tau $$ 
$$ \le  \int_0^s\Big( \big| 2\,\sh\!^2\big[{\ts\frac{w_\tau}{2}}\big] - {\ts\frac{w^2_\tau}{2}} \big| + \big| 2\,\sh\!^2\big[{\ts\frac{\beta_\tau}{2}}\big] - {\ts\frac{\beta^2_\tau}{2}} \big| + \O\big(R^2\tau^{2}\big) \Big) d\tau  $$ 
$$ \le  \int_0^s\Big( R^4 \tau^{2} + \O\big(R^2\tau^{2}\big)\Big) d\tau\, \le \,  R^4 s^{3} \,\le \, s^{5/2} \quad \hbox{ for }\, 0\le s\le s_R'' >0\, , $$ 
with probability $\, 1-\O(e^{- {R^{2/3}/4}})\,$. \   In particular we can take $\,\kappa = 2/9\,$ in the statement.  $\;\diamond$   

\parm \parmn 
\ub{\bf  Proof of Lemma \ref{lem.LaplZ}} \quad  Let us use \big([Y], Chapter (2) ``The laws of some quadratic functionals of Brownian motion''\big),  
considering for any $\,b> 0\,$ the exponential martingale 
\begin{equation*}
M^b_s := \exp \!\Big(- {\ts\frac{b}{2}}\, \big(w_s^2-w_0^2-s\big) - {\ts\frac{b^2}{2}}\int_0^s w_\tau^2\, d\tau\Big) = \exp \!\Big(- b\int_0^sw_\tau\,dw_\tau - {\ts\frac{b^2}{2}}\int_0^s w_\tau^2\, d\tau\Big), 
\end{equation*}
and the new probability $\,\P^b$ having on $\FF_s\,$ density $M_s^b\,$ with respect to $\P$.  As noticed in \big([Y], (2.1.1)\big), Girsanov's Theorem yields a real $(\P^b,\FF_s)$-Brownian motion $(B_u\,,\,0\le u\le s)$ such that \  ${\ds w_u = w_0 + B_u - b\int_0^u w_\tau\, d\tau\,}$, which means that under $\,\P^b$,  $\,w\,$ has become an Ornstein-Uhlenbeck process, alternatively expressed by \   ${\ds w_u = e^{-b\,u} \bigg( w_0+ \int_0^u e^{b\,\tau}\, dB_\tau\bigg) }$.  Therefore 
\begin{equation*}
\E_0\bigg[ \exp\!\bigg(\rt1\! \int_0^s (a + c\, \tau)\, dw_\tau -  {\ts\frac{b^2}{2}} \int_0^s w_\tau^2\, \, d\tau\bigg)\bigg] 
\end{equation*}
\begin{equation*}
= \, e^{-b\,s/2}\times \E^b_0\bigg[ \exp\!\bigg(\rt1\! \int_0^s (a + c\, \tau)\, dw_\tau + b\, w_s^2/2\bigg)\bigg]
\end{equation*}
on the one hand, and on the other hand for any test-function $\,f\,$ on $\R\,$: 
\begin{equation*}
\int_0^s f(u)\, dw_u = \int_0^s f(u)\, \bigg[ dB_u - b\, e^{-b\,u} \bigg(\int_0^u e^{b\,\tau}\, dB_\tau \bigg) du  \bigg]  
\end{equation*}
\begin{equation*}
= \int_0^s \bigg( f(\tau ) - b\, e^{b\,\tau}\! \int_\tau^s f(u)\, e^{-b\,u} du \bigg) dB_\tau \, ,   
\end{equation*}
so that 
\begin{equation*}
\E^b_0\Bigg[ \bigg(\int_0^s f(\tau)\, dw_\tau\bigg)^2\Bigg] = \int_0^s \Big[ f(\tau ) - b\, e^{b\,\tau}\! \int_\tau^s f(u)\, e^{-b\,u} du \Big]^2\, d\tau \, . 
\end{equation*}
Taking \  $f(\tau) = a + c\, \tau\,$, we obtain 
\begin{equation*}
\E^b_0\Bigg[ \bigg(a\, w_s + c\int_0^s \tau\, dw_\tau\bigg)^2\Bigg] = \int_0^s \Big[ a + c\, \tau + (a+c/b)(e^{b (\tau-s)}-1) + c\, (s\, e^{b (\tau-s)}-\tau) \Big]^2\, d\tau 
\end{equation*}
\begin{equation*}
= b\2 \int_0^s \Big[ (ab+c + bc\, s)\, e^{b(\tau-s)} - c \Big]^2\, d\tau 
\end{equation*}
\begin{equation*}
= b\2 \int_0^s \Big[ (ab+c + bc\, s)^2\, e^{2b (\tau-s)} - 2c\,(ab+c + bc\, s)\, e^{b (\tau-s)} +c^2 \Big] d\tau 
\end{equation*}
\begin{equation*}
= b^{-3} \Big[ \5 \big(b\, a+ (bs+1)c\big)^2\, (1-e^{-2b s}) - 2c\,\big(b\, a+ (bs+1)c\big)\,(1-e^{-b s}) + bs\, c^2 \Big]   
\end{equation*}
\begin{equation*}  {\ts 
=\,  { \frac{1-e^{-2b s}}{2\, b}}\,a^2 + { \frac{bs -1+ 2 e^{-b\, s} - (bs+1) e^{-2b s}}{b^2}}\, ac
+ { \frac{b^2s^2-3 + 4\,(bs+1)\, e^{-b s} - (bs+1)^2\, e^{-2b s}}{2\, b^3}}\, c^2\, . }
\end{equation*}
This yields the covariance matrix of the $\,\P^b_0$-Gaussian variable ${\ds \Big( w_s \, , \int_0^s \tau\, dw_\tau\Big)}$, hence its joint law. Namely the covariance matrix under $\,\P^b_0\,$ of ${\ds \Big( \sqrt{2b}\, w_s \, , \, \sqrt{2b^3} \int_0^s \tau\, dw_\tau\Big)}$ is 
\begin{equation*}
K_{bs} = \begin{pmatrix} 1-e^{-2b s}\; & bs -1+ 2 e^{- b s} - (bs+1) e^{-2b s}  \cr 
bs -1+ 2 e^{- b s} - (bs+1) e^{-2b s}\; & b^2s^2-3 + 4\,(bs+1)\, e^{-bs} - (bs+1)^2\, e^{-2bs} 
\end{pmatrix}\!  , 
\end{equation*} 
and its determinant is \  $\delta_{bs} := 2(bs-2) + 8\, e^{-b\,s} - 2(bs+2)\, e^{-2b\,s}$, which increases with $\,bs>0\,$ and then does not vanish.  
\if{  and its inverse is 
\begin{equation*}
K_{bs}\1 = \delta_{bs}\1 \begin{pmatrix} b^2s^2-3 + 4\,(bs+1)\, e^{-b s} - (bs+1)^2\, e^{-2b s} & 
1 - bs - 2 e^{- b s} + (bs+1) e^{- 2 b s}  \cr 
1 - bs - 2 e^{- b s} + (bs+1) e^{- 2 b s}  & 1-e^{-2b s} \end{pmatrix}\!  , 
\end{equation*} 
}\fi  
Therefore the density of ${\ds \Big( \sqrt{2b}\, w_s \, , \, \sqrt{2b^3} \int_0^s \tau\, dw_\tau\Big)}$ is 
\begin{equation*}
(u,v)\lmt \, \frac{1}{2\pi\sqrt{\delta_{bs}}}\, \exp\bigg[ \frac{-1}{2\,\delta_{bs}} \Big(  \alpha_{bs}\, u^2 + 2\,\gamma_{bs} \, uv +  \big(1-e^{-2b s}\big)\, v^2\Big) \bigg] ,
\end{equation*}
with \   $\,\alpha_x := x^2-3 + 4\,(x+1)\, e^{-x} - (x+1)^2\, e^{-2x}\,$ and \  $ \gamma_x := 1 - x - 2 e^{- x} + (x+1) e^{- 2 x}$, \parn 
so that \  $\delta_x = (1-e^{-2x})\, \alpha_x - \gamma_x^2\,$.  \   Hence  
\begin{equation*}
\E^b_0\bigg[ \exp\!\bigg(\rt1\! \int_0^s (a + c\, \tau)\, dw_\tau + b\, w_s^2/2\bigg)\bigg]
\end{equation*}
\begin{equation*}
= \frac{1}{2\pi\sqrt{\delta_{bs}}} \int_{\sR^2} e^{ \frac{a\rt1}{\sqrt{2b}}\, u + \frac{c\rt1}{\sqrt{2b^3}}\, v + u^2/4} \exp\bigg[ \frac{-1}{2\,\delta_{bs}} \Big(  \alpha_{bs}\, u^2 + 2\,\gamma_{bs} \, uv +  \big(1-e^{-2b s}\big)\, v^2\Big) \bigg]  du\, dv 
\end{equation*}
\begin{equation*}
= \frac{1}{\sqrt{2\pi (1-e^{-2b s})}} \int_{\sR} e^{ \frac{a\rt1}{\sqrt{2b}}\, u + \frac{\delta_{bs}}{2(1-e^{-2b s})} \big(\frac{c\rt1}{\sqrt{2b^3}} - \frac{\gamma_{bs}}{\delta_{bs}}\, u\big)^2+ \big( \frac{1}{4}- \frac{\alpha_{bs}}{2\,\delta_{bs}}\big) u^2 } \,  du
\end{equation*}
\begin{equation*}
= \frac{e^{\frac{- \delta_{bs}\, c^2}{4(1-e^{-2b s}) b^3}}}{\sqrt{2\pi (1-e^{-2b s})}} 
\int_{\sR} e^{ \big(\frac{a\rt1}{\sqrt{2b}} - \frac{c\rt1 \gamma_{bs}}{(1-e^{-2b s})\sqrt{2b^3}}\big)\, u - \big( \frac{\alpha_{bs}}{\delta_{bs}} - \frac{1}{2} - \frac{\gamma_{bs}^2}{(1-e^{-2b s})\delta_{bs}} \big) u^2/2 } \,  du
\end{equation*}
\begin{equation*}
= \frac{e^{\frac{- \delta_{bs}\, c^2}{4(1-e^{-2b s}) b^3}}}{\sqrt{2\pi (1-e^{-2b s})}} 
\int_{\sR} e^{ \big(a - \frac{c\, \gamma_{bs}}{b\,(1-e^{-2b s})}\big)\frac{\rt1}{\sqrt{2b}}\,  u - \big( \frac{1+e^{-2b s}}{1-e^{-2b s}}\big) u^2/4 } \,  du
\end{equation*}

\begin{equation*}
= \frac{2\, e^{\frac{- \delta_{bs}\, c^2}{4(1-e^{-2b s}) b^3}}}{\sqrt{2 (1+e^{-2b s})}} 
\times 
e^{\frac{-1}{2b^3}\times \frac{1-e^{-2b s}}{1+e^{-2b s}} \big(a b - \frac{c\, \gamma_{bs}}{1-e^{-2b s}}\big)^2} 
\end{equation*}
\begin{equation*}
= \frac{e^{bs/2}}{\sqrt{\ch\!(b s)}} \, \exp\!\bigg[\frac{-1}{4b^3}\bigg( \frac{\delta_{bs}\, c^2}{1-e^{-2b s}} + 2\, \frac{1-e^{-2b s}}{1+e^{-2b s}} \Big(a\, b - \frac{\gamma_{bs}\, c}{1-e^{-2b s}}\Big)^2\bigg)\bigg] 
\end{equation*}
\begin{equation*}
= \frac{e^{bs/2}}{\sqrt{\ch\!(b s)}} \, \exp\!\bigg[\frac{-1}{4b^3 (1+e^{-2b s})}\Big(  2\, (1-e^{-2b s})\, b^2 a^2 + (2\alpha_{bs}-\delta_{bs})\, c^2 -  4\,b\, \gamma_{bs}\, ac\Big)\bigg] . 
\end{equation*}
Therefore for any $\,b, s > 0$ and real $a,c\,$ we obtain\,: 
\begin{equation*}
\E_0\bigg[ \exp\!\bigg(\rt1\! \int_0^s (a + c\, \tau)\, dw_\tau -  {\ts\frac{b^2}{2}} \int_0^s w_\tau^2\, \, d\tau\bigg)\bigg] 
\end{equation*}
\begin{equation*}
= \frac{1}{\sqrt{\ch\!(b s)}} \, \exp\!\bigg[\frac{-1}{4\,b^3 (1+e^{-2b s})}\Big(  2\, (1-e^{-2b s})\, b^2 a^2 + (2\alpha_{bs}-\delta_{bs})\, c^2 -  4\, \gamma_{bs}\,b\, ac\Big)\bigg] . 
\end{equation*}
Then taking \  $a=r+c\,s \,$, we get$\,$: 
\begin{equation*}
\E_0\bigg[ \exp\!\bigg(\rt1\Big[ r\, w_s + c \int_0^s w_\tau\, d\tau \Big]-  {\ts\frac{b^2}{2}} \int_0^s w_\tau^2\, \, d\tau\bigg)\bigg] 
\end{equation*}
\begin{equation*}
=\, \E_0\bigg[ \exp\!\bigg(\rt1\! \int_0^s (r + c\, s - c\,\tau)\, dw_\tau -  {\ts\frac{b^2}{2}} \int_0^s w_\tau^2\, \, d\tau\bigg)\bigg] 
\end{equation*}
\begin{equation*}
= \frac{1}{\sqrt{\ch\!(b s)}} \, \exp\!\bigg[- \frac{ 2 (1-e^{-2b s})\, b^2 (r+c\,s)^2 + (2\alpha_{bs}-\delta_{bs})\, c^2 + 4\, \gamma_{bs}\,b\, (r+c\,s)\, c}{4\,b^3 (1+e^{-2b s})}\bigg]  
\end{equation*}
\begin{equation*}
=\,  \frac{1}{\sqrt{\ch\!(b s)}} \exp\!\bigg[- \frac{ (1-e^{-2b s})\, b^2 (r^2 + 2s\, r c) +  \big({bs}\,(1+e^{-2b s})-(1-e^{-2b s})\big) c^2 + 2b\, \gamma_{bs}\, r c}{2\,b^3 (1+e^{-2b s})}\bigg]  
\end{equation*}
\begin{equation*}
=\,  \frac{1}{\sqrt{\ch\!(b s)}} \, \exp\!\bigg[- \frac{\th\!(b s)}{2\, b} \,r^2 - \frac{ {bs} -\th\!(b s)}{2\,b^3}\, c^2
- 2\, \frac{ \sh\!^2(b s/2)}{b^2\, \ch\!(b s)}\, r c\bigg] . \;\; \diamond 
\end{equation*}

\pars \parmn
\ub{\bf Proof of Proposition \ref{pro.abscCv}} \quad Denote by $\,\lambda_0\in\R_+\,$ the abscissa of convergence of the integral, so that the map \    ${\ds \lambda\mapsto \int_{0}^\infty  e^{\lambda\, x} \; q_1(w,\beta , x , \zeta , z) \, dx}\;$ is analytic on $\big\{\Re(\lambda)<\lambda_0\big\}$. By Proposition \ref{pro.LaplY} and Lemma \ref{lem.explPsi} it is equal to $\,\Phi(\lambda)$ for $\,\Re(\lambda) <  \min\{4\pi^2,\lambda_0\}$. \   
Hence, for any real $\,\lambda<\min\{4\pi^2,\lambda_0\}$ and $\,t\in\R\,$ we have \vspace{-3mm} 
\begin{equation*}
\Phi(\lambda+\rt1 t) = \int_{0}^\infty  e^{\rt1 t\, x}\, e^{\lambda\, x} \; q_1(w,\beta , x , \zeta , z) \, dx\, .  \end{equation*}  
Let us show now that \  $t\mapsto \Phi(\lambda+\rt1 t)$ belongs to $L^1\cap L^2(\R)$, in order to inverse the above Fourier transform. \  Of course we have to deal here with the large values of $\,|t|\,$, i.e., of $\,|\lambda+\rt1 t|\,$ in the expression (\ref{f.Psi-lambd})\,:  \vspace{-3mm} 
\begin{equation*}  
\Phi_{w,\beta,\zeta,z}(\lambda) =\, \frac{\lambda^2\, e^{B' \lambda - \frac{\sqrt{\lambda}}{4}\,\cotg\!({\frac{\sqrt{\lambda}}{2}})(w^2+\beta^2)}}{8\pi^2 \big[1- \cos({\sqrt{\lambda}}\,) - (\sqrt{\lambda}/2) \sin({\sqrt{\lambda}}\,)\big]} \, \exp\!\Bigg[\frac{B^2 \lambda}{1- \frac{2}{\sqrt{\lambda}}\,\tg\!(\frac{\sqrt{\lambda}}{2}) }\Bigg] ,
\end{equation*}
in which $\,\lambda\,$ is to be replaced by $\,\lambda+\rt1 t$, and we have set \ $\,B^2:= \frac{(w\!-\!2z)^2\! + \! (\beta\!-\!2\zeta)^2}{8}\,$ and $\,B':= \frac{4wz + 4\beta\zeta -w^2-\beta^2}{8} = \frac{z^2+\zeta^2}{2} - B^2$\raise0.2pt\hbox{.} \    Then for any  $\,t\in\R\,$ we have\,:  \vspace{-2mm} 
\begin{equation*}  {\ts 
\frac{\sqrt{\lambda+\rt1 t}}{2}\, =: \alpha + \rt1 b = \sqrt{\frac{|\lambda+\rt1 t| +\lambda}{8}} + \rt1 {\rm sign}(t) \sqrt{\frac{|\lambda+\rt1 t| -\lambda}{8}} }  \vspace{-2mm} 
\end{equation*}
and  \vspace{-1mm}  
\begin{equation*}
\tg\!\Big[{\ts\frac{\sqrt{\lambda+\rt1 t}}{2}}\Big] =\frac{\tg \alpha + \rt1 \th b}{1- \rt1\! \tg \alpha\, \th b} = \frac{\sin(2\alpha) + \rt1\! \sh\!(2b)}{\ch\!(2b) + \cos(2\alpha)}\,= \rt1\! {\rm sign}(t) +\O\Big(e^{-\sqrt{|t|/2}}\,\Big) .  
\end{equation*}
Therefore, for large $|t|$ we have$\,$:  \vspace{-2mm}
\begin{equation*} 
\exp\!\left[B'(\lambda+\rt1 t) - (w^2+\beta^2)\, {\ts\frac{\sqrt{\lambda+\rt1 t}}{4}\,\cotg\!\big({\frac{\sqrt{\lambda+\rt1 t}}{2}}\big)}+ {\ts \frac{B\, (\lambda+\rt1 t)}{1- \frac{2}{\sqrt{\lambda+\rt1 t}}\,\tg\!\big(\frac{\sqrt{\lambda+\rt1 t}}{2}\big)}}\right] 
\end{equation*}
\begin{equation*}
=\, e^{\big(\frac{z^2+\zeta^2}{2}\big) \lambda} \exp\!\bigg[{\big({\ts\frac{z^2+\zeta^2}{2}}\big) \rt1\! t + \rt1 {\rm sign}(t)({\ts\frac{w^2+\beta^2}{4}}+2B)\sqrt{\lambda\! +\!\rt1\! t} \Big[ 1\! + \O\big(|t|^{-1/2}\big)\Big]}\! \bigg] ,
\end{equation*}
the modulus of which is  \vspace{-2mm} 
\begin{equation*}
\, e^{({z^2+\zeta^2})\, \lambda/2}\, \exp\!\Big[{- ({\ts\frac{w^2+\beta^2}{4}}+2B^2) \sqrt{\ts\frac{|\lambda+\rt1 t| -\lambda}{2}}\, + \O(1)\Big]} . 
\end{equation*}
Moreover  
\begin{equation*}
\cos\!\big({\ts{\sqrt{\lambda+\rt1 t}}}\,\big) = \ch\!(2b)\cos(2\alpha) - \rt1 \sh\!(2b)\sin(2\alpha) 
\end{equation*}
and
\begin{equation*}
\sin\!\big({\ts{\sqrt{\lambda+\rt1 t}}}\,\big) = \ch\!(2b)\sin(2\alpha) + \rt1 \sh\!(2b)\cos(2\alpha)
\end{equation*}
entail  
\begin{equation*}
\Big|1- \cos\!\big({\ts{\sqrt{\lambda+\rt1 t}}}\,\big) - \big({\ts\sqrt{\lambda+\rt1 t}}/2\big) \sin\!\big({\ts{\sqrt{\lambda+\rt1 t}}}\,\big)\Big|^2 
\end{equation*}
\begin{equation*}
= [\ch\!(2b)-\cos(2\alpha)]^2+ (\alpha^2+b^2)[\ch\!^2(2b)-\cos^2(2\alpha)] - 2[\ch\!(2b)-\cos(2\alpha)][b\,\sh\!(2b)+\alpha\sin(2\alpha)]
\end{equation*} 
\begin{equation*}  
= [\ch\!(2b)-\cos(2\alpha) - b\,\sh\!(2b) - \alpha\sin(2\alpha)]^2+ [\alpha\,\sh\!(2b)-b\sin(2\alpha)]^2
\end{equation*}
\begin{equation*}
= (\alpha^2+b^2)\, \sh\!^2(2b) + \O\big( \alpha\,\sh\!^2(2b)\big) = \Big[{\ts\frac{|\lambda+\rt1 t|}{4}} + \O\big(\sqrt{|t|}\,\big)\Big]\, \sh\!^2\sqrt{\ts\frac{|\lambda+\rt1 t| -\lambda}{2}}\, \quad \hbox{for large $\,|t|\,$}\raise0pt\hbox{.} 
\end{equation*}
So far, for $\,\lambda< 4\pi^2$ and for large $\,|t|\,$ we have\,: \    $\big|\Phi_{w,\beta,\zeta,z}(\lambda+\rt1 t)\big| $ 
\begin{equation*}
=\, (\lambda^2+t^2)^{3/4}\, e^{({z^2+\zeta^2})\, \lambda/2}\, \exp\!\Big[{- ({\ts\frac{w^2+\beta^2}{4}}+2B^2+1) \sqrt{\ts\frac{\sqrt{\lambda^2+t^2} -\lambda}{2}}\,+ \O(1)\Big]} \Big[1+ \O\big(|t|^{-1/2}\big)\Big]  
\end{equation*}
\begin{equation*}
=\, \O\Big(e^{2\pi^2({z^2+\zeta^2})}\Big) |t|^{3/2} \exp\!\Big[- \big({\ts\frac{w^2+\beta^2}{4}}+2B^2+1+ \O(|t|\1)\big) \sqrt{|t|/2}\,\Big] 
\end{equation*}
\begin{equation*}
=\, \O\Big(e^{2\pi^2({z^2+\zeta^2})}\Big) \exp\!\Big[- \big({\ts\frac{w^2+\beta^2}{4}}+2B^2+1/2\big) \sqrt{|t|/2}\,\Big] .   
\end{equation*}
Hence we can inverse the above Fourier transform for $\,\lambda<\min\{4\pi^2,\lambda_0\}$, and thus we obtain the wanted (\ref{f.invFourq1}), which holds a posteriori for $\,\lambda< 4\pi^2$\,:  \vspace{-2mm} 
\begin{equation*}    
e^{\lambda\, x} \, q_1(w,\beta , x , \zeta , z) = \int_{-\infty}^\infty  e^{-\rt1 t\, x}\, \Phi(\lambda+\rt1 t) \, \frac{dt}{2\pi}\, \raise1.7pt\hbox{,} \; \hbox{ for } \, x>0  \, \hbox{ and for} \; \lambda< 4\pi^2 . 
\end{equation*}
Thence, for any real $\,\lambda< 4\pi^2$, taking $\,\e< 4\pi^2-\lambda\,$ for positive $\,x\,$ we have  \vspace{-2mm} 
\begin{equation*}
e^{\lambda\, x} \, q_1(w,\beta , x , \zeta , z) = \O\big(e^{-\e\, x}\big)\times  \int_{-\infty}^\infty \big| \Phi(\lambda+\e+\rt1 t)\big| \, dt\, = \O\big(e^{-\e\, x}\big) ,
\end{equation*}
which entails the integrability of $\,x\mapsto e^{\lambda\, x} \, q_1(w,\beta , x , \zeta , z)$, so that  finally $\,\lambda_0\ge 4\pi^2$. $\diamond$  \par

\if{                   

\section{Expected result about the density of the process $X_s$} \label{sec.Resultat?} \indf 
     In analogy with the Li-Yau-like Gaussian bounds obtained in [DM], where the ``tangent process'' is Gaussian, we should get here bounds where the Gaussian density is replaced by the density of the tangent process $\,Y_s\,$. Actually, we expect to have a result as follows. 
\bthe \label{th.secchbound} \   The law of the $\R^5$- valued Dudley diffusion $(X_s)$, which satisfies  (in the coordinates specified in Section \ref{sec.RepHR2}) the system of stochastic differential equations (\ref{f.sde})(\ref{f.sde'}), admits a density  $\,\tilde p_s(X) = \tilde p_s(\lambda\,,\, \mu\,,\, x\,,\, y\,,\, z)$ \big(when started from $X_0\equiv 0$, and with respect to the Lebesgue measure $\,\Lambda(dX) = d\lambda\, d\mu \, dx\,dy\,dz$\big) such that for any $S>0$, there exists a constant $\,C=C_S>1$ such that for any $\,s\in\, ]0,S]$ and any $\,X=(\lambda,\mu,x,y,z)\in\R^5\,$: 
$$ \frac{1}{C\, s^{6}}\;  q_1^{C}\bigg( \frac{\lambda}{\sqrt{s}}\, \raise1.8pt\hbox{,}\, \frac{\mu}{\sqrt{s}}\, \raise1.8pt\hbox{,}\,  \frac{x-s}{s^2}\, \raise1.8pt\hbox{,}\, \frac{y}{\sqrt{s^3}}\, \raise1.8pt\hbox{,}\, \frac{z}{\sqrt{s^3}} \bigg) \le \, \tilde p_s(X) \le\,  \frac{C}{s^{6}}\;  q_1^{\frac{1}{C}}\bigg( \frac{\lambda}{\sqrt{s}}\, \raise1.8pt\hbox{,}\, \frac{\mu}{\sqrt{s}}\, \raise1.8pt\hbox{,}\,  \frac{x-s}{s^2}\, \raise1.8pt\hbox{,}\, \frac{y}{\sqrt{s^3}}\, \raise1.8pt\hbox{,}\, \frac{z}{\sqrt{s^3}} \bigg) ,
$$ 
where the density $\,q_1$ of $\,Y_1\,$ is given by (\ref{f.p1=}). 
\ethe

\parn 
Le 11 oct. 2011, Shizan Fang a \'ecrit : \quad Cher Jacques,

J'ai lu un peu ton livre pendant l'\'et\'e: il semble que la diffusion horizontale sur $\mathrm{PSO}(1,d)$ v\'erifie m\^eme la condition hypoelliptique de Bismut (les crochets d'ordre 1 suffisent \`a engendrer les espaces tangents) ; elle permet de construire la diffusion de Dudley dans le chapitre suivant.\parsn
Ma r\'eponse : \par 

   La diffusion de Dudley est construite dans la section VII.6.7, et est \'equivalente \`a la diffusion 
matricielle $(X_s)$ sur le groupe de Poincar\'e, solution de (VII.12) ; \parn
les $(A_j\,|\,j=1,\ldots,d)$ et leurs crochets d'ordre 1 engendrent  lin\'eairement $\,\mathrm{so}(1,d)$, et d'autre part $\,A_0 \equiv e_0\,$ et les $[A_j,A_0] \equiv e_j$ (pour $j=1,\ldots,d$) engendrent  $\R^{1,d}$ ; ce qui fait le compte, en effet$\,$: les $(A_j\,|\,j=1,\ldots,d)$ et leurs crochets d'ordre 1 suffisent bien \`a engendrer lin\'eairement l'alg\`ebre de Lie $\,\GG\,$ du groupe de Poincar\'e. 

En revanche sur $T^1\H^d$ directement, et donc dans un autre syst\`eme de coordonn\'ees, 
il faut un crochet d'ordre 2 pour engendrer l'espace tangent (ceci pour tout $d>1$). 
Cela semble indiquer que cette fois encore il vaut peut-être mieux rester au niveau du 
fibr\'e des rep\`eres, c'est-\`a-dire du groupe de Poincar\'e, sans projeter. 

Je ne sais toutefois pas si c'est un gain d\'efinitivement appr\'eciable, car la difficult\'e 
semble provenir pour une bonne part du fait que dans les formules du genre Campbell-Hausdorff
les crochets it\'er\'es (il y en a peu qui sont nuls) produisent de toute façon une infinit\'e de termes ; 
nous ne sommes pas du tout dans un cas r\'esoluble (puisque $[\GG,\GG] = \GG$).  

Encore une pr\'ecision :  \   pour avoir une densit\'e  $p_s(...)$, il faut l'hypoellipticit\'e de $(\partial_s - L)$, i.e. il faut pouvoir se passer de $\{A_0\}$ isol\'e, pour engendrer $\GG$ ;  
et l\`a un crochet d'ordre 2 est indispensable, justement pour r\'ecup\'erer $A_0\,$. 
De sorte que la condition de Bismut me semble insuffisante$\,$;  non ?

}\fi                  

\vskip 8mm

\centerline{\bfseries \LARGE References} \label{Ref} 
\par \vskip 6mm 


\vbox{ \noindent 
{\bf [A]} \   Azencott R.  \quad {\it Densité des diffusions en temps petit : d\'eveloppements asymptotiques.} \par \hskip 30mm   Sém. Proba. XVI (1980-81), Lecture Notes n$^o\,$921, 237-284, Springer. }
\par \bigskip 

\vbox{ \noindent 
{\bf [BA1]} \  Ben Arous G.  \quad {\it Développement asymptotique du noyau de la chaleur hypoelliptique \par \hskip 32mm   hors du cut-locus.}  \  Ann. sci. \'E.N.S., Sér. 4, t. 21 n$^o\,$3, 307-331, 1988. }
\par \bigskip 

\vbox{ \noindent 
{\bf [BA2]} \  Ben Arous G.  \quad {\it Développement asymptotique du noyau de la chaleur hypoelliptique \par \hskip 38mm   sur la diagonale.}  \  Ann. Inst. Fourier 39, n$^o\,$1, 73-99, 1989. }
\par \bigskip 

\vbox{ \noindent 
{\bf [BF]} \  Bailleul I, Franchi J.   \ {\it  Non-explosion criteria for relativistic diffusions.}  \par \hskip 46mm  
To appear at Annals of Probability. }
\par \bigskip 



\vbox{ \noindent 
{\bf [BPY]} \  Biane P., Pitman J., Yor M. \quad {\it Probability laws related to the Jacobi theta and \par \hskip 9mm    Riemann zeta functions, and Brownian excursions.} \par  \hskip 10mm Bull.  A.M.S. vol. 38, n$^o\,$4, 435-465, 2001. }
\par \bigskip 


\vbox{ \noindent 
{\bf [Ca]} \  Castell F. \ \  {\it Asymptotic expansion of stochastic flows.}
\par \smallskip \hskip 27mm  Prob. Th. Rel. Fields 96, 225-239, 1993. }
\par \bigskip 

\vbox{ \noindent 
{\bf [Co]} \  Copson E.T. \ \  {\it Asymptotic Expansions.}
\hskip 7mm  Cambridge University Press, 1965. }
\par \bigskip 

\vbox{ \noindent 
{\bf [CDJR]} \  Chan T., Dean D.S.,  Jansons K.M., Rogers L.C.G.   \ {\it  On Polymer Conformations  \par \hskip 33mm  in Elongational Flows. }  \hskip 1mm  Commun. Math. Phys. 160, 239-257, 1994. }
\par \bigskip 

\vbox{ \noindent 
{\bf [DM]} \  Delarue F., Menozzi S. \ \  {\it  Density estimates for a random noise propagating through\par \hskip 33mm  a chain  of differential equations. }  \hskip 5mm   J. F. A.   $n^o\,$259, 1577-1630, 2010. }
\par \bigskip 

\vbox{ \noindent 
{\bf [Du]} \  Dudley R.M. \ \  {\it Lorentz-invariant Markov processes in relativistic phase space.}
\par \smallskip \hskip 33mm  Arkiv f\"or Matematik 6, n$^o\,$14, 241-268, 1965. }
\par \bigskip 


\vbox{ \noindent 
{\bf [ERS]} \  ter Elst A.F.M., Robinson D.W., Sikora A. \  {\it Small time asymptotics of diffusion \par  \hskip 40mm processes.} \quad J. Evol. Equ. 7, n$^o\,$1, 79-112, 2007. }
\par \bigskip 

\vbox{ \noindent 
{\bf [F-LJ1]} \  Franchi J., Le Jan Y.  \ \  {\it Relativistic Diffusions and Schwarzschild Geometry.} \par
\hskip 54mm Comm. Pure Appl. Math., vol. LX, n$^o\,$2, 187-251, 2007.}
\par \bigskip 

\vbox{ \noindent 
{\bf [F-LJ2]} \   Franchi J., Le Jan Y.  \ \  {\it Curvature diffusions in general relativity.} \par
\hskip 55mm   Comm. Math. Physics, vol. 307, n$^o\,$2, 351-382, 2011. } 
\par \bigskip 

\vbox{ \noindent 
{\bf [F-LJ3]} \   Franchi J., Le Jan Y.  \ \   {\it Hyperbolic Dynamics and Brownian Motion}. \par \hskip 12mm   Oxford Mathematical Monographs,  Oxford Science Publications, august 2012.  } 
\par \bigskip 




%
\vbox{ \noindent 
{\bf [L]} \   Léandre R. \quad {\it  Intégration dans la fibre associée à une diffusion dégénérée.} \par \hskip 28mm    Prob. Th. Rel. Fields  76,  n$^o\,$3, 341-358, 1987.}
\par \bigskip 

\vbox{ \noindent 
{\bf [NSW]} \  Nagel A., Stein E.M., Wainger S.  \   {\it  Balls and metrics defined by vector fields I\,: \par \hskip 33mm   Basic properties.}  \qquad    Acta Math. 155, 103-147, 1985. }
\par \bigskip 

%
\vbox{ \noindent 
{\bf [T]} \  Tardif C. \quad {\it A Poincar\'e cone condition in the Poincar\'e group.}
\par \smallskip \hskip 25mm {\footnotesize To appear at} Potential Analysis. }
\par \bigskip 

\vbox{ \noindent 
{\bf [V]} \  Varadhan S.R.S. \ \  {\it  Diffusion Processes in a Small Time Interval.}
\par \hskip 37mm  Comm. Pure Applied Math.   n$^o\,$20, 659-685, 1967. }
\par \bigskip 

%
\vbox{ \noindent 
{\bf [Y]} \  Yor M. \quad {\it Some aspects of Brownian motion. Part I. Some special functionals.}
\par \smallskip \hskip 21mm Lectures in Mathematics, ETH Zürich. Birkhäuser Verlag, Basel, 1992. }
\par \bigskip 

\end{document}